\pdfoutput=1
\documentclass[11pt, letterpaper]{amsart}
\makeatletter

\usepackage[T1]{fontenc}
\usepackage{lmodern}
\usepackage{textcmds}  
\usepackage{amsmath, amssymb, amsfonts, amstext, verbatim, amsthm, mathrsfs, stmaryrd}
\usepackage{microtype}
\usepackage[all,cmtip,2cell]{xy}
\usepackage{enumerate}
\usepackage{enumitem}
\usepackage[colorlinks=true,linkcolor=blue,citecolor=blue,urlcolor=blue,citebordercolor={0 0 1},urlbordercolor={0 0 1},linkbordercolor={0 0 1}]{hyperref} 
\usepackage[shortalphabetic]{amsrefs} 
\usepackage[nameinlink]{cleveref}
\usepackage{tikz-cd}
\usepackage{enotez}
\setenotez{backref=true}
\usepackage[normalem]{ulem}

 \usepackage{pbox}

 \usepackage[colorinlistoftodos]{todonotes}

\tikzset{
  symbol/.style={
    draw=none,
    every to/.append style={
      edge node={node [sloped, allow upside down, auto=false]{$#1$}}}
  }
}

\def\makeCal#1{%
\expandafter\newcommand\csname c#1\endcsname{\mathcal{#1}}}
\def\makeBB#1{%
\expandafter\newcommand\csname b#1\endcsname{\mathbb{#1}}}
\def\makeFrak#1{%
\expandafter\newcommand\csname f#1\endcsname{\mathfrak{#1}}}

\count@=0
\loop
\advance\count@ 1
\edef\y{\@Alph\count@}%
\expandafter\makeCal\y
\expandafter\makeBB\y
\expandafter\makeFrak\y
\ifnum\count@<26
\repeat

\theoremstyle{plain}
\newtheorem{thm}{Theorem}[section]
\newtheorem{cor}[thm]{Corollary}
\newtheorem{lem}[thm]{Lemma}
\newtheorem{prop}[thm]{Proposition}

\theoremstyle{definition}
\newtheorem{rem}[thm]{Remark}
\newtheorem{defn}[thm]{Definition}
 \newtheorem{notn}[thm]{Notation}
 \newtheorem{ex}[thm]{Example}

\newtheorem{introthm}{Theorem}

\DeclareMathOperator{\Spec}{Spec}
\DeclareMathOperator{\Aut}{Aut}

\DeclareMathOperator{\op}{op}
\DeclareMathOperator{\Branch}{Branch}
\DeclareMathOperator{\Coh}{Coh}
\DeclareMathOperator{\Pic}{Pic}
\DeclareMathOperator{\GL}{GL}
\DeclareMathOperator{\Alg}{Alg}
\DeclareMathOperator{\Chow}{Chow}

\DeclareMathOperator{\reduced}{\red}
\DeclareMathOperator{\wt}{wt}

\DeclareMathOperator{\Hom}{Hom}

\DeclareMathOperator{\Quot}{Quot}
\DeclareMathOperator{\Tot}{Tot}

\DeclareMathOperator{\Sym}{Sym}

\DeclareMathOperator{\Span}{Span}
\DeclareMathOperator{\Proj}{Proj}

\DeclareMathOperator{\pr}{pr}
\DeclareMathOperator{\disc}{disc}

\DeclareMathOperator{\Cond}{Cond}
\DeclareMathOperator{\Hilb}{Hilb}

\newcommand{\LS}{\mathcal{G}}
\DeclareMathOperator{\NS}{NS}
\DeclareMathOperator{\QCoh}{QCoh}
\DeclareMathOperator{\Perf}{Perf}
\DeclareMathOperator{\Map}{Map}

\newcommand{\CM}{L_{\mathrm{CM}}}

\newcommand{\basefield}{k}
\newcommand{\red}{\mathrm{red}}

\newcommand{\univ}{\mathrm{univ}}

\makeatletter

\title{Projectivity of the moduli of equidimensional branchvarieties}
\author{Daniel Halpern-Leistner, Andres Fernandez Herrero,\\ Trevor Jones, and Ritvik Ramkumar}

\address{(Daniel Halpern-Leistner) Department of Mathematics\\Cornell University\\Ithaca, NY\\USA}
\email{daniel.hl@cornell.edu}
\address{(Andres Fernandez Herrero) Department of Mathematics\\University of Pennsylvania\\Philadelphia, PA\\USA}
\email{andresfh@sas.upenn.edu}
\address{(Trevor Jones) Department of Mathematics\\ Johns Hopkins University\\Baltimore, Maryland\\USA}
\email{ttj9@cornell.edu}
\address{(Ritvik Ramkumar) Department of Mathematics\\University of Notre Dame\\South Bend, IN\\USA}
\email{rramkuma@nd.edu}

\begin{document}

\begin{abstract}
    We resolve an open problem posed by Alexeev-Knutson on the projectivity of the moduli of branchvarieties in the equidimensional case. As an application, we construct projective moduli spaces of reduced equidimensional varieties equipped with ample linear series and subject to a semistability condition.
\end{abstract}

\maketitle
\tableofcontents

\section{Introduction}
The classification of (reduced) varieties has been a central problem in algebraic geometry for more than a century. A modern avatar of such classification questions is the construction of moduli spaces of varieties, which has been an active and fruitful endeavor over the last decades. Mumford's work on Geometric Invariant Theory (GIT) \cite{mumford1994geometric} provided a powerful tool for the construction of such moduli spaces. This was illustrated by Gieseker in his GIT construction of the moduli of curves \cite{gieseker_moduli_curves}, as well as in more recent generalizations, such as the moduli of stable maps \cite{stble_maps_git}. GIT constructions also feature prominently in the construction of moduli of elliptic surfaces \cite{miranda-elliptic-p1, miranda-elliptic-lectures}.

The main GIT approach to constructing moduli spaces of varieties, initiated by Mumford in his lecture notes \cite{mumford-stability-projective-varieties}, is to consider Grothendieck's Hilbert scheme $\Hilb^P(\bP^N)$ \cite[Part IV]{fga_grothendieck}, which is a projective scheme parametrizing flat families of subschemes of $\mathbb{P}^N$ with Hilbert polynomial $P$. Then one uses GIT to construct the moduli space as an orbit space for the action of $\mathrm{SL}_N$ on a locally closed subscheme of the Hilbert scheme. However, it has proved challenging to connect GIT semistability on Hilbert schemes with an intrinsic notion of semistability, i.e., one that does not depend on the choice of embedding. As a result, the field has turned to different methods for constructing moduli spaces in the special case of general type \cite{kollar-families-varieties-book} and log-Fano varieties \cite{liu-xu-zhuang, xu-fano-book}.

In \cite{branchvarieties-paper}, Alexeev and Knutson introduced the stack of branchvarieties, $\Branch^P(X)$, as an alternative to the Hilbert scheme. $\Branch^P(X)$ is a proper Deligne-Mumford stack that parametrizes flat families of projective, reduced, equidimensional schemes with Hilbert polynomial $P$ equipped with a finite morphism to a fixed projective scheme $X$. In \cite{branchvarieties-paper}, a branchvariety is not required to be equidimensional, but in this paper we only address the equidimensional case.

Unlike the Hilbert scheme, the stack $\Branch^P(\bP^N)$ compactifies the moduli of reduced subschemes of $\bP^N$ without introducing non-reduced and non-equidimensional schemes. It would seem that $\Branch^P(\mathbb{P}^N)$, or its coarse moduli space, is better suited for the construction of certain moduli of varieties using GIT. However, in their 2010 paper \cite{branchvarieties-paper}, Alexeev and Knutson left as an open question whether the coarse moduli space of $\Branch^P(X)$ is projective, which is a necessary hypothesis to set up a GIT problem.

The main goal in this paper is to prove this conjecture, thus opening up a new avenue for GIT constructions of moduli spaces of varieties.

\begin{introthm}\label{T:projectivity_branch}
     Let $k$ be a field of characteristic $0$, and let $X$ be a projective $k$-scheme equipped with an ample line bundle $\cO_X(1)$. Fix a Hilbert polynomial $P$ of degree $n$. Then the moduli stack $\Branch^P(
    X)$ of equidimensional branchvarieties with Hilbert polynomial $P$ has a projective coarse moduli space.
\end{introthm}

\begin{rem}
    We do not assume that the field $k$ is algebraically closed for any of the statements in this introduction. However, note that ultimately all of our arguments reduce easily to the case when $k$ is algebraically closed. Indeed, our proof proceeds by showing that certain line bundles defined over the ground field $k$ are (semi)ample, and we may check this by passing to the algebraic closure of $k$. 
\end{rem}

\subsection*{Moduli of ample linear series}

As an application of our main theorem, we use GIT to construct moduli spaces of ample linear series. Specifically, we describe a moduli stack $\LS_P^r$ that parameterizes reduced equidimensional polarized schemes $(X,\cO_X(1))$ with fixed Hilbert polynomial $P(t)$ along with a surjection $\cO_X \otimes V \to \cO_X(1)$, where $V$ is a vector space of dimension $r+1$. Using the theory of $\Theta$-stability, we introduce a notion of polynomial semistability in \Cref{D:polynomial_stability}, analogous to $K$-semistability of varieties, which implies that $V \to \Gamma(X,\cO_X(1))$ is injective. We prove the following:

\begin{introthm}[= \Cref{T:moduli_ample_linear_series}]
    The stack $\LS_P^{r,\rm{ss}}$ of semistable ample linear series admits a projective good moduli space $M_P^{r}$.
\end{introthm}

The moduli spaces $M_P^r$ appear to be closely related to moduli spaces of limit linear series on curves, which have had many applications to Brill-Noether theory and the geometry of the moduli of curves \cite{eisenbud-harris, osserman-linear-series}. It would be interesting to investigate this relationship further.

\subsection*{Line bundles on Branch and Hilb.}
\Cref{T:projectivity_branch} is a consequence of a detailed analysis of ampleness of certain line bundles on $\Branch^P(X)$.

For any Noetherian scheme $T$ and any $T$-flat family of proper schemes $\pi: Y \to T$ equipped with a morphism $f: Y \to X$, we define the following sequence of line bundles indexed by an integer $r \in \mathbb{Z}$
\[ \lambda_r := \det R\pi_*(f^*(\cO_X(r))).\]
The formation of the line bundles $\lambda_r$ is compatible with base-change, and therefore it defines corresponding line bundles on $\Hilb^P(X)$ and $\Branch^P(X)$, which we also denote by $\lambda_r$. By the work of Mumford and Knudsen \cite{knudsen-mumford}, it is known that the assignment $r \mapsto \lambda_r$ is a polynomial with values in the Picard group of $\Branch^P(X)$. More precisely, for $0 \leq i \leq n+1$, there are line bundles $b_i$ on $\Branch^P(X)$ such that we have 
\[ \lambda_r = \sum_{i=0}^{n+1} b_i \binom{r}{i}\]
for all $r \in \mathbb{Z}$, where we use additive notation for addition in the Picard group.\\

\subsection*{Chow morphisms.}
Suppose that $\cO_X(1)$ is a very ample line bundle on $X$ inducing a closed immersion $X \subset \mathbb{P}^N$. The assignment given by pushing forward the fundamental class induces well-defined morphisms from $\Branch^P(X)$ and $\Hilb^P(X)$ to an appropriate Chow scheme of cycles on $\mathbb{P}^N$. This induces morphisms $\Chow: \Branch^P(X) \to \bP(V_{P,N})$ and $\Chow: \Hilb^P(X) \to \bP(V_{P,N})$, where $V_{P,N}$ is a $k$-vector space that depends only on $P$ and $N$. See \Cref{subsection: chow morphism} for a more detailed discussion. In \Cref{prop: b_{n+1} is the pullback of O(1) under Chow}, we observe that $\Chow^\ast(\cO(1)) \cong b_{n+1}$, thus implying that $b_{n+1}$ is semiample. Our main strategy will be to show that this Chow morphism is relatively projective.

\begin{introthm}[= \Cref{prop: -lambda_n is relatively ample} + \Cref{prop: -b_n is relatively semiample} + \Cref{prop: semiampleness of a b_n+1 - b_n}]\label{introthm: ample_line_bundles}
     Fix a Hilbert polynomial $P$ of degree $n$. Let $X$ be a projective $k$-scheme equipped with a very ample line bundle $\cO_X(1)$ which induces an embedding $X \hookrightarrow \mathbb{P}^N$. Then the following hold:
\begin{enumerate}
    \item For every sufficiently large positive integer $r\gg0$, there exists $a(r)>0$ such that for all rational numbers $a\geq a(r)$ the $\mathbb{Q}$-line bundle $-a b_n -\lambda_r$ is relatively ample for the Chow morphism $\Chow: \Branch^P(X) \to \bP(V_{P,N})$.
    \item The line bundle $-b_n$ is relatively semiample for the Chow morphism $\Chow: \Branch^P(X) \to \bP(V_{P,N})$.
    \item For any rational number $a>n$, the $\mathbb{Q}$-line bundle $ab_{n+1} - b_n$ is semiample on $\Branch^P(X)$.
    \end{enumerate}
    \end{introthm}
    
A surprising feature of \Cref{introthm: ample_line_bundles} is the fact that $ab_n + \lambda_r$ is relatively \emph{anti}ample for $\Chow: \Branch^P(X) \to \bP(V_{P,N})$ for $r \gg 0$. This is in contrast with the case of $\Hilb^P(X)$. Indeed, it is known that $\lambda_r$ is ample on $\Hilb^P(X)$ for $r \gg0$, and therefore, it is relatively ample for $\Chow: \Hilb^P(X) \to \bP(V_{P,N})$.\\

\subsection*{Refinements over the deminormal locus.} In addition, we prove the following stronger ampleness statement on the open substack $\Branch_{\mathrm{Dn}}^P(X) \subset \Branch^P(X)$ parametrizing deminormal branchvarieties (i.e. equidimensional branchvarieties which satisfy Serre's condition $S_2$ and have nodal singularities in codimension $1$).
    
\begin{introthm}[=\Cref{thm: ampleness of cm line bundle on s_2 locus}] \label{introthm: S_2 locus}
        Fix a Hilbert polynomial $P$ of degree $n$. Let $X$ be a projective $k$-scheme equipped with a very ample line bundle $\cO_X(1)$. Let $a>n$ be a rational number. Let $f: \Branch^P(X) \to \mathbb{P}^M$ denote any morphism induced by a base-point free linear system of a power of the semiample $\mathbb{Q}$-line bundle $ab_{n+1} - b_n$. Then the restriction of $f$ to $\Branch_{\mathrm{Dn}}^P(X) \subset \Branch^P(X)$ is quasi-finite. In particular, the $\mathbb{Q}$-line bundle $ab_{n+1} - b_n$ is ample on $\Branch_{\mathrm{Dn}}^P(X)$.
    \end{introthm}

\begin{ex}
    A one-dimensional scheme is deminormal if and only if it has nodal singularities. Hence, when $\deg(P)=1$ and $a>1$, the line bundle $a b_{n+1} - b_n$ is ample on the open stack $\Branch^P_{\mathrm{nod}}(X)$ of branchvarieties whose source is a nodal curve.
\end{ex}

\subsection{Notation and conventions} \label{notn: conventions}
In this paper, we work over a fixed ground field $k$ of characteristic $0$. Unless otherwise stated, all stacks and schemes are defined over $k$. For any two stacks $\cX$ and $\cY$, an undecorated product $\cX \times \cY$ will always denote a fiber product over $k$. We may sometimes write $\cX_{\cY}$ to denote $\cX \times \cY$.

\begin{defn}
    We say that a morphism $f: \cX \to Y$ from an algebraic stack $\cX$ to a Noetherian scheme $Y$ is \textbf{projective} if $f$ is proper with finite relative inertia, and the relative coarse space $X \to Y$ of $\cX$ is projective over $Y$.
\end{defn} 

\begin{defn}
    If $f: \cX \to Y$ is a projective morphism and $M$ is a rational line bundle in $\Pic(\cX)\otimes \bQ$, then we say that $M$ is \textbf{relatively ample} (resp. \textbf{relatively semiample}) if, locally on $Y$, a sufficiently divisible power descends to a relatively ample (resp. relatively semiample) line bundle on $X \to Y$.
\end{defn}

If $Y = \Spec(k)$, then we say that the stack $\cX$ is projective and the rational line bundle $M$ is ample (resp. semiample) on $\cX$.

\subsection{Acknowledgements}

The first named author was supported by NSF grants DMS-1945478 and DMS-2052936, and a Simons Foundation fellowship. He would like to thank the Simons Laufer Mathematical Sciences Institute for its hospitality while this research was conducted. The fourth named author was partially supported by NSF grant DMS-2401462. We would like to thank Jarod Alper, Allen Knutson, Yuji Odaka, David Rydh, and Chenyang Xu for helpful discussions while this work was in progress. We would also like to thank an anonymous referee for thoughtful suggestions.



\section{The stack of branchvarieties}

\subsection{Pure algebras}

\begin{defn}[{\cite[Defn. 2.2]{halpern2024moduli}, \cite[Defn. 1.1.2]{huybrechts2010geometry}}] \label{defn: absolute pure} 
Let $Y$ be a scheme of finite type over a field and let $y \in Y$ be a point. We say that $y$ has dimension $n$ if the closure $\overline{y} \subset Y$ is a scheme of dimension $n$. 
Let $F$ be a coherent $\mathcal{O}_Y$-module. We say that $F$ is a \textbf{pure sheaf} of dimension $n$ if all the associated points of $F$ have dimension $n$.
\end{defn}

\begin{defn}{\cite[Defn. 2.3]{halpern2024moduli}} \label{defn: relative pure} 
Let $\pi: Y \to T$ be a finite type morphism of schemes. We say that a sheaf $F$ on $Y$ is \textbf{$T$-pure} of dimension $n$ if it is $T$-flat, locally finitely presented as a $\cO_T$-module, and for all $t \in T$ we have that $F_t$ is a pure sheaf of dimension $n$ on $Y_{t}$.
\end{defn}

\begin{defn} \label{defn: pure coherent algebra}
Let $\pi: Y \to T$ be a finite-type morphism of schemes. A $T$-pure \textbf{$\cO_Y$-algebra} $(\cA, m, u)$ of dimension $n$ is the data of
\begin{enumerate}
    \item[(i)] a $T$-pure sheaf $\cA$ on $Y$ of dimension $n$,
    \item[(ii)] a multiplication morphism $m : \cA^{\otimes 2} \to \cA$, and
    \item[(iii)] a unit morphism $u : \cO_Y \to \cA$
\end{enumerate}
that endow $\cA$ with the structure of a $\cO_Y$-sheaf of commutative unital algebras.
A $T$-pure $\cO_Y$-algebra $\cA$ is said to be reduced if for all geometric points $t \in T$ the scheme $\Spec_{Y_t}(\cA|_{Y_t})$ is reduced.
\end{defn}

\begin{defn}[Moduli of pure algebras] 
Let $X$ be a projective $\basefield$-scheme.
We define $\Alg^{n}(X)$ to be the pseudofunctor from $(\text{Sch}/\basefield)^{\op}$ to groupoids that associates to every $\basefield$-scheme $T$ the groupoid of $T$-pure $\mathcal{O}_{X_T}$-algebra of dimension $n$.
\end{defn}

\begin{prop} \label{prop: stack of pure algebras is algebraic}
Let $X$ be a projective $\basefield$-scheme. The pseudofunctor $\Alg^{n}(X)$ is an algebraic stack with affine diagonal and locally of finite type over $\basefield$.
\end{prop}

\begin{proof}

There is a forgetful morphism $\Alg^{n}(X) \to \Coh^n(X)$, where $\Coh^n(X)$ is the stack parameterizing $n$-dimensional families of pure sheaves \cite[Defn. 2.4]{halpern2024moduli}. Since the stack $\Coh^{n}(X)$ is algebraic, with affine diagonal and locally of finite type over $\basefield$ \cite[Prop. 2.5]{halpern2024moduli}, it suffices to show that the forgetful morphism is affine and of finite type.

Choose a scheme $T$ and a morphism $T \to \Coh^n(X)$ corresponding to a family of pure sheaves $\cA$ on $X_T$. A $T'$ point of the fiber product $\Alg^{n}(X) \times_{\Coh^n(X)} T$ is the data of a commutative unital algebra structure on the base change $\cA_T$, which consists of a section $u: \cO_{X_{T'}} \to \cA_{T'}$ and a morphism $m: \cA_{T'}^{\otimes 2} \to \cA_{T'}$. Since $\cA$ is a $T$-flat and $\cA$ and $\cA^{\otimes 2}$ are finitely presented $\cO_{X_T}$-modules, by \cite[Lem. 2.16]{halpern2024moduli} the functor 
\[ H := \Hom_{T}(\cA^{\otimes 2}, \cA) \times_T \Hom_{T}(\cO_T, \cA) \] is representable by a scheme which is affine and of finite presentation over $T$. Furthermore, the subfunctor of $H$ parametrizing sections $(u,m)$ which satisfy the axioms of a commutative unital algebra structure is representable by a closed subscheme of $H$.
\footnote{Indeed, the algebra axioms amount to the equality of certain pairs of universally-defined morphisms of sheaves on $X_H$ (the morphisms lie in $\Hom(A_H, A_H)$ for the unit condition, in $\Hom(A^{\otimes 2}_H, A_H)$ for the commutativity condition, and in $\Hom(A^{\otimes 3}_H, A_H)$ for associativity). These morphisms correspond to pairs of sections of the affine morphisms $\Hom_H(A_H, A_H) \to H$, $\Hom_{H}(A^{\otimes 2}_H, A_H) \to H$ and $\Hom_{H}(A^{\otimes 3}_H, A_H) \to H$ respectively \cite[Lem. 2.16]{halpern2024moduli}. Since these morphisms are separated, the subfunctor of $H$ parametrizing points where these pairs of sections agree is represented by a closed immersion.}

 We conclude that $T \times_{\Coh^n(X)} \Alg^n(X) = Z$ is represented by a relatively affine scheme of finite type over $T$, as desired.
\end{proof}

\begin{notn}
    Fix an ample line bundle $\cO_X(1)$ on $X$. For any given polynomial $P$, we denote by $\Alg^P(X)$ the stack whose $T$-points are $T$-pure $\cO_{X_T}$-algebras $\cA$ such that $\cA_t$ viewed as a sheaf on $X_t$ has Hilbert polynomial $P$ for all $t \in T$. 
\end{notn}

Since the Hilbert polynomial is locally constant in flat families of sheaves, $\Alg^{n}(X)$ can be written as a disjoint union of open and closed substacks $\Alg^{n}(X) = \bigsqcup_{\deg(P) = n} \Alg^{P}(X)$.

\begin{defn}
    We denote by $\Alg^{n}(X)_{\reduced}$ the substack of $\Alg^n(X)$ parametrizing families of reduced algebras as in \Cref{defn: pure coherent algebra}. The inclusion $\Alg^{n}(X)_{\reduced}\hookrightarrow \Alg^n(X)$ is an open immersion by {\cite[\href{https://stacks.math.columbia.edu/tag/0C0E}{Tag 0C0E}]{stacks-project}}.
\end{defn}

\begin{defn}[Determinant line bundles] \label{defn: determinantal line bundle}
Fix an ample line bundle $\cO_X(1)$ on $X$.
Let $\cA_{\univ}$ denote the universal sheaf on $\Alg^{n}(X) \times X$, and let $\pi: \Alg^n(X) \times X \to \Alg^n(X)$ 
denote the first projection. 
For every integer $r \in \mathbb{Z}$, we set 
$\lambda_r := \text{det} \, R\pi_{\ast}\left(\cA_{\univ}(r) \right) \in \Pic(\Alg^n(X))$.
\end{defn}

The morphism $\Spec_{\Alg^n(X) \times X}(\cA_{\univ}) \to \Alg^n(X)$ is flat and proper, and the condition of a complex on a stack being perfect is local. Therefore, $R\pi_\ast(\cA_{\univ})$ is a perfect complex on $\Alg^n(X)$ and its formation commutes with base change \cite[\href{https://stacks.math.columbia.edu/tag/0A1P}{Tag 0A1P}]{stacks-project}. The determinant of perfect complexes on schemes was constructed in \cite[Thm.~2, pp.~42]{knudsen-mumford}, and the construction was extended to perfect complexes on stacks in \cite[Sect.~3.1]{schurg2015derived}. In fact, the latter constructs a morphism of higher stacks $\det : \underline{\Perf} \to B\bG_m$, where $\underline{\Perf}$ is the stack that takes any scheme $T$ to the $\infty$-groupoid of perfect complexes on $T$ and quasi-isomorphisms between them. In particular, one gets a morphism of groupoids
\[\xymatrix@R=5pt{
\Perf(\Alg^n(X)) \ar@{=}[d] \ar[r]^{\det} & \Pic(\Alg^n(X)) \ar@{=}[d] \\ \Map(\Alg^n(X), \underline{\Perf}) \ar[r] & \Map(\Alg^n(X),B\bG_m).
}\]
The fact that $\underline{\Perf} \to B\bG_m$ is a morphism of stacks implies that the formation of $\lambda_r$ is compatible with base change.

\begin{defn}[Mumford-Knudsen coefficients] \label{defn: knudsen-mumford coefficients}
By \cite[Thm. 4]{knudsen-mumford}, the assignment $r \mapsto \lambda_r$ is a polynomial in the variable $r$ of degree $n+1$ with values in $\Pic(\Alg^n(X))$. More precisely,
$$
\lambda_r = \sum_{i = 0}^{n+1} \binom{r}{i}b_i 
    = \frac{r^{n+1}}{(n+1)!}b_{n+1} 
     + \left({\frac{1}{n!}}b_n - \frac{1}{2(n-1)!}b_{n+1} \right)r^n
     +  O(r^{n-1})
$$
for certain line bundles $b_i$. 
\end{defn}

\begin{rem}[Formula for $b_i$] \label{remark: taylor series expansion}
    Using iterated difference operators  \cite[Pg. 54-55]{knudsen-mumford}, the line bundles in \Cref{defn: knudsen-mumford coefficients} can be expressed as 
$$ 
b_i = \sum_{j = 0}^i {(-1)^{i-j} \binom{i}{j}}\lambda_{j}.
$$ 
\end{rem}

\begin{rem}
    The statement in \cite[Thm. 4]{knudsen-mumford} is stated in the case when the base is a Noetherian scheme, whereas our base $\Alg^n(X)$ is an algebraic stack. However, we note that a direct application of smooth descent shows that the equality in \Cref{defn: knudsen-mumford coefficients} holds with $b_i$ defined as in \Cref{remark: taylor series expansion}. Indeed, we can apply the result in \cite[Thm. 4]{knudsen-mumford} after passing to an atlas $\sqcup_i U_i \to \Alg^n(X)$ where each $U_i$ is a Noetherian scheme, since we are working with families of sheaves supported in dimension $n$, which satisfy condition $Q_{n}$ from \cite[pg. 50]{knudsen-mumford}. Since the isomorphisms obtained in \cite[Thm. 4]{knudsen-mumford} are canonical and functorial, the isomorphisms over the atlas $\sqcup_i U_i$ automatically satisfy the cocycle condition required for smooth descent. Hence they glue to yield the desired identifications as in \Cref{defn: knudsen-mumford coefficients}.
\end{rem}

\begin{prop} \label{prop: det and b_n} Let $p: \Spec(\basefield) \to \bP^n$ be a $k$-point. Let $\det(\cA_{univ}) \in \Pic(\Alg^n(\mathbb{P}^n) \times \mathbb{P}^n)$ denote the determinant of the universal sheaf, as in e.g. \cite[pg. 36,37]{huybrechts2010geometry}.
The line bundles $p_{\Alg^n(\bP^n)}^{\ast}(\det(\cA_{\emph{univ}}))$ and $b_n$ on $\Alg^n(\bP^n)$ are isomorphic.
\end{prop}
\begin{proof}
    Let $\cO_p$ be the skyscraper sheaf of $\mathbb{P}^n$ at the closed point $p$. Then $\cO_p \cong \cO_p \otimes \cO(n)$ has a Koszul resolution whose $i^{th}$ term is $\cO(n) \otimes \bigwedge^i (\cO(-1)^{\oplus n}) \cong \cO(n-i) \otimes_k \bigwedge^i(k^{\oplus n})$. It follows from the formula in \Cref{remark: taylor series expansion} that $b_n \cong \det (R \pi_\ast(\tau^*\cO_p \otimes^L \cA_{\univ}))$, where $\tau: \Alg^n(X) \times X \to X$ is the second projection and $\otimes^L$ is the derived tensor product. On the other hand, the functoriality of $\det(-)$ under derived pullback of perfect complexes implies that $p_{\Alg^n(X)}^\ast(\det(\cA_{\univ})) \cong \det (Lp_{\Alg^n(X)}^\ast(\cA_{\univ}))$. The lemma therefore follows from the claim that the canonical morphism $L p_{\Alg^n(X)}^\ast(\cA_{\univ}) \to R \pi_\ast(\tau^*\cO_p \otimes^L \cA_{\univ})$ is a quasi-isomorphism. This last claim can be checked after applying $(p_{\Alg^n(X)})_\ast$, at which point it follows from the quasi-isomorphism $\tau^*\cO_p \otimes^L \cA_{\univ} \cong (p_{\Alg^n(X)})_\ast(Lp_{\Alg^n(X)}^\ast(\cA_{\univ}))$.
    \footnote{Replace $\tau^*\cO_p \otimes^L \cA_{\univ}$ with $(p_{\Alg^n(X)})_\ast(Lp^\ast(\cA_{\univ}))$, then observe that $$(p_{\Alg^n(X)})_\ast \pi_\ast (p_{\Alg^n(X)})_\ast Lp_{\Alg^n(X)}^\ast(\cA_{\univ}) \cong (p_{\Alg^n(X)})_\ast Lp_{\Alg^n(X)}^\ast(\cA_{\univ})$$ because $p_{\Alg^n(X)} \circ \pi \circ p_{\Alg^n(X)} = p_{\Alg^n(X)}$.}
\end{proof}

\subsection{The discriminant morphism}
\begin{notn}[{\cite[\href{https://stacks.math.columbia.edu/tag/0BVH}{Tag 0BVH}]{stacks-project}}] \label{notn: discriminant}
    Given a Noetherian scheme $S$ and a locally free sheaf of commutative unital $\cO_S$-algebras $\cA$ on $S$, there is a symmetric bilinear trace pairing $\cA \otimes \cA \to \cO_S$ defined by sending a pair of sections $a,b \in \cA(U)$ over an open $U \subset S$ to the trace of the morphism of locally free sheaves $(ab) \cdot (-): \cA|_{U} \to \cA|_{U}$ induced by multiplication by the product $ab$. We may interpret the trace pairing as a morphism $\cA \to \cA^{\vee}$, which yields a morphism of line bundles $\det(\cA) \to \det(\cA)^{\vee}$. Tensoring by $\det(\cA)^{\vee}$, we a canonically defined morphism $s: \cO_S \to \det(\cA)^{-2}$ called the discriminant section. The vanishing locus $D_{\cA} := V(s) \hookrightarrow S$ is called the discriminant of $\cA$. The formation of $s$ is compatible with arbitrary base-change on $S$, and the complement of $S \setminus D_{\cA}$ is the maximal open subscheme over which the morphism $\Spec_{S}(\cA) \to S$ is \'etale.
\end{notn}

More generally, let $T$ be a Noetherian scheme over $k$, and choose a $T$-point of $\Alg^n(\mathbb{P}^n)_{\reduced}$ corresponding to a $T$-flat family of torsion-free sheaves on $\mathbb{P}^n_T$ equipped with a commutative unital algebra structure. Let $U_{\cA} \subset \mathbb{P}^n_T$ denote maximal open subscheme where $\cA$ is locally free. For every point $t \in T$, the $t$-fiber of the closed complement $\left(\mathbb{P}^n_T \setminus U_{\cA}\right)_t$ has codimension at least $2$ \cite[Lemma. 2.5]{rho-sheaves}. Let $\det(\cA)$ denote the determinant line bundle of the family $\cA$ \cite[pg. 36,37]{huybrechts2010geometry}. By \Cref{notn: discriminant}, we have a canonical section $s: \cO_{U_{\cA}} \to \det(\cA)^{-2}|_{U_{\cA}}$ defined over the open $U_{\cA}$. By Hartogs's theorem \cite[Lemma 2.3(c)]{rho-sheaves}, this extends uniquely to a canonical section $\widetilde{s}: \cO_{\mathbb{P}_T^n} \to \det(\cA)^{-2}$ whose formation commutes with arbitrary base-change on $T$.
\begin{defn}[Discriminant divisor] \label{defn: discriminant divisor}
Given a $T$-point $\cA: T \to \Alg^n(\mathbb{P}^n)_{\reduced}$ as above, the \textbf{discriminant} $D_{\cA} \subset \mathbb{P}^n_T$ of $\cA$ is defined to be the vanishing locus $V(\widetilde{s}) \subset \mathbb{P}^n_T$ of the unique extension $\widetilde{s}: \cO_{\mathbb{P}^n_T} \to \det(\cA)^{-2}$.
\end{defn}

\begin{rem} \label{remark: etale locus in terms of discriminant locus}
    It follows from \Cref{notn: discriminant} that $U_{\cA} \cap (\mathbb{P}^n_T \setminus D_{\cA})$ is the largest open subscheme of $\mathbb{P}^n_T$ over which the morphism $\Spec_{\mathbb{P}^n_T}(\cA) \to \mathbb{P}^n_T$ is \'etale. 
\end{rem}

\begin{prop} \label{prop: discriminant morphism algred}
    Fix a degree $n$ Hilbert polynomial $P(t) = \sum_{i=0}^n a_i \binom{t}{i}$ and set $d := 2na_n-2a_{n-1}$. The assignment $\cA \mapsto D_{\cA}$ induces a well-defined morphism $\disc: \Alg^P(\mathbb{P}^n)_{\reduced} \to \mathbb{P}\left(H^0(\mathbb{P}^n, \cO(d)\right)$ to the projective space of degree $d$ relative effective Cartier divisors on $\mathbb{P}^n$. Furthermore, we have $\disc^*(\cO(1)) = -2b_n$.
\end{prop}
\begin{proof}
    For any Noetherian scheme $T$ and any morphism $\varphi: T \to \Alg^P(\mathbb{P}^n)_{\reduced}$ corresponding to a $T$-family of reduced algebras $\cA$ on $\mathbb{P}^n_T$, consider the discriminant $D_{\cA}$ cut out by a section $s: \cO_{\mathbb{P}^n_T} \to \det(\cA)^{-2}$. The line bundle $\det(\cA)$ on $\mathbb{P}^n_T$ is isomorphic to $\cO(m) \otimes \pi^*(L)$, where $m \in \mathbb{Z}$ is an integer, $\pi: \mathbb{P}^n_T \to T$ denotes the structure morphism, and $L \in \Pic(T)$. 
     Since the Hilbert polynomial of $\cA$ is $P(t) = \sum_{i=0}^n a_i \binom{t}{i}$, it follows that $m = a_{n-1} - na_n$.
    \footnote{The integer $m$ is the degree of $\cA$. If we write the Hilbert polynomial in the standard polynomial basis, i.e., 
     $
     P(t) = \sum_{i=0}^na'_i\frac{t^i}{i!},
     $
     then $\deg(\cA) = a'_{n-1} - \frac{a'_n(n+1)}{2}$ \cite[Definition 1.2.11]{huybrechts2010geometry}.
     If we instead use the binomial basis, it follows that $a'_n = a_n$ and 
     $
     a'_{n-1} = -\frac{a_{n}(n-1)}{2} + a_{n-1}. 
     $
     Thus,
     $$
     m = \deg \cA = -\frac{a_{n}(n-1)}{2}+a_{n-1} - \frac{a_n(n+1)}{2} 
         = a_{n-1} -na_n. 
     $$}
    Therefore, $D_{\cA}$ is cut out by a section of $\cO(d) \otimes \pi^*(L)^{-2}$. For every $t \in T$, the scheme $\Spec_{\mathbb{P}^n_t}(\cA_t)$ is reduced. Since the characteristic of the ground field $k$ is $0$, it follows that the finite morphism $\Spec_{\mathbb{P}^n_t}(\cA_t) \to \mathbb{P}^n_t$ is \'etale over a nonempty open subscheme of $\mathbb{P}^n_t$. By the fiber-wise criterion for flatness \cite[\href{https://stacks.math.columbia.edu/tag/05VJ}{Tag 05VJ}]{stacks-project}, it follows that $\Spec_{\mathbb{P}^n_T}(\cA) \to \mathbb{P}^n_T$ is flat over an open subscheme $U \subset \mathbb{P}^n_T$ which intersects every $T$-fiber, and after shrinking $U$ we have that $\Spec_{\mathbb{P}^n_T}(\cA) \to \mathbb{P}^n_T$ is \'etale over $U$. By \Cref{remark: etale locus in terms of discriminant locus}, it follows that the Cartier divisor $D_{\cA} \subset \mathbb{P}^n_T$ does not contain any $T$-fiber, and therefore it is a relative effective Cartier divisor \cite[\href{https://stacks.math.columbia.edu/tag/062Y}{Tag 062Y}]{stacks-project}. The formation of $D_{\cA}$ commutes with base-change, and hence the assignment $\cA \mapsto D_{\cA}$ induces a well-defined morphism $\disc: \Alg^P(\mathbb{P}^n)_{\reduced} \to \mathbb{P}\left( H^0(\mathbb{P}^n, \cO(d) )\right)$. Furthermore, since $D_{\cA}$ is cut out by a section of $\cO(d) \otimes \pi^*(L)^{-2}$, it follows that $\varphi^* \disc^*(\cO(1)) = L^{-2}$. 
    Hence to conclude the proof it suffices to show that there is a canonical identification $L = \varphi^*(b_n)$. The choice of a section $p: \Spec(k) \to \mathbb{P}^n$ yields canonical identifications $L = p^*(\det(\cA)) = \varphi^*(b_n)$, where the right-most isomorphism was constructed in \Cref{prop: det and b_n}.
\end{proof}

\subsection{Moduli of equidimensional branchvarieties}



\begin{defn}
    Let $X$ be a projective $k$-scheme, and let $T$ be any $k$-scheme. A \textbf{$T$-family of equidimensional branchvarieties} of dimension $n$ over $X$ is a finite morphism $f:Y \to X_T$ such that $Y \to T$ is flat of finite presentation and every fiber $f_t: Y_t \to X_t$ is reduced and equidimensional.
\end{defn}

\begin{defn} \label{defn: branch varieties}
Let $X$ be a projective $k$-scheme. The stack $\Branch^n(X)$ of equidimensional branchvarieties of dimension $n$ is the pseudofunctor from $(\text{Sch}/k)^{\text{op}}$ to groupoids that associates to every $k$-scheme $T$ the groupoid of $T$-families of equidimensional branchvarieties of dimension $n$ over $X$.
\end{defn}

\begin{defn}
    If $X$ is equipped with an ample line bundle $\cO_X(1)$, then for any given polynomial $P(t)$ we denote by $\Branch^P(X) \subset \Branch^n(X)$ the open and closed subfunctor parametrizing families of branchvarieties such that $(X_t,\cO_{X_t}(1))$ has Hilbert polynomial $P$ for any $t\in T$.
\end{defn}

The main theorem of \cite{branchvarieties-paper} contains the following as a special case.
\begin{thm}[{\cite[Theorem 0.4]{branchvarieties-paper}}] \label{thm: Branch is proper}
    Let $X$ be a projective $k$-scheme equipped with a very ample line bundle $\cO_X(1)$. Fix a Hilbert polynomial $P$. The stack $\Branch^P(X)$ is a proper Deligne-Mumford stack. In particular, it admits a proper coarse moduli space.
\end{thm}

\begin{lem} \label{lemma: branchn isomorphic to algred}
Let $X$ be a projective $k$-scheme equipped with a polarization $\cO_X(1)$. There is an isomorphism of stacks $\psi: \Alg^{n}(X)_{\reduced} \xrightarrow{\sim} \Branch^{n}(X)$ which restricts to an isomorphism $\psi: \Alg^P(X)_{\reduced} \xrightarrow{\sim} \Branch^P(X)$ for any given Hilbert polynomial $P$.
\end{lem}
\begin{proof} There is a natural equivalence at the level of $T$-points. 
Indeed, given a family of reduced algebras $A$ on $X_T$, we obtain a family $\Spec_{X_T}(A) \to X_T$ of branchvarieties over $T$. Conversely, given a family of branchvarieties $f : Y \to X_T$ over $T$, the $\cO_{X_T}$-algebra $f_\ast(\cO_Y)$ is a family of reduced algebras. 
This is evidently compatible with the corresponding Hilbert polynomials.
\end{proof}

\begin{notn}
    We shall abuse notation and also write $\lambda_r \in \Pic(\Branch^n(X))$  for the pullback $(\psi^{-1})^*(\lambda_r)$ of the line bundle $\lambda_r$ from \Cref{defn: determinantal line bundle}. Note that this is compatible with the notation in \cite[\S5]{branchvarieties-paper}. Similarly, we denote by $b_i \in \Pic(\Branch^n(X))$  the pullback of the line bundle $b_i$ from \Cref{defn: knudsen-mumford coefficients}. 
\end{notn}

\begin{cor} \label{coroll: discriminant morphism branch}
    Fix a degree $n$ Hilbert polynomial $P(t) = \sum_{i=0}^n a_i \binom{t}{i}$ and set $d := 2na_n-2a_{n-1}$. There is a discriminant morphism $\disc: \Branch^P(\mathbb{P}^n) \to \mathbb{P}\left(H^0(\mathbb{P}^n, \cO(d))\right)$ such that $\disc^*(\cO(1)) = -2b_n$.
\end{cor}
\begin{proof}
    This is a consequence of \Cref{lemma: branchn isomorphic to algred} and \Cref{prop: discriminant morphism algred}.
\end{proof}

\begin{defn}
    Let $Z \subset X$ be a closed subscheme of a projective $k$-scheme $X$. We denote by $\Branch^n(X \setminus Z)$ the open substack of $\Branch^n(X)$ parametrizing $T$-families of branchvarieties $f: Y \to X_T$ such that the image of $f$ lies on the open subscheme $X_T \setminus Z_T$.
\end{defn}

\begin{lem} \label{lemma: linear projection}
    Let $\mathbb{P}^{N - n-1} \subset \mathbb{P}^N$ be a linear subspace. Then there is a well-defined linear projection morphism $\pr: \Branch^n(\bP^{N} \setminus \bP^{N - n - 1}) \to \Branch^n(\mathbb{P}^n)$ which is affine of finite type and satisfies $\pr^*(\lambda_r) = \lambda_r$.
\end{lem}
\begin{proof}
Fix a linear subspace $\bP^n \subset \mathbb{P}^N$ disjoint from the given $\bP^{N-n-1}$. Let $f: Y \to \mathbb{P}^N_T \setminus \mathbb{P}^{N-n-1}_T$ be a family of branchvarieties parametrized by a scheme $T$. The linear projection $\pr: \mathbb{P}^N_T \setminus \mathbb{P}^{N-n-1}_T \cong \Tot_{\bP^n_T}(\cO_{\bP^n_T}(1)^{\oplus (N-n)}) \to \bP^n_T$ is an affine morphism. The composition $\pr \circ f$ is affine, since $f$ is finite and $\pr$ is affine. It is also proper, since $Y$ and $\bP^n_T$ are proper over $T$. It follows that $\pr \circ f:Y \to \mathbb{P}^n_T $ is finite and, thus, a family of equidimensional branchvarieties. Therefore, the assignment $f \mapsto \pr \circ f$ defines a morphism $\pr:\Branch^n(\bP^{N} \setminus \bP^{N - n - 1}) \to \Branch^n(\mathbb{P}^n)$.

Since it is evident that $f^*(\cO_{\bP_T^N}(1)) = (\pr \circ f)^*(\cO_{\bP_T^n}(1))$,
it follows that
$$
\lambda_r = \det R\pi_{\ast}f^*(\cO_{\bP^N}(r)) = \det R\pi_{\ast}(\pr \circ f)^*(\cO_{\bP^n}(r)) = \pr^{\ast}(\lambda_r).
$$

We are only left to show that $\pr: \Branch^n(\bP^{N} \setminus \bP^{N - n - 1}) \to \Branch^n(\mathbb{P}^n)$ is affine and of finite type. Choose a scheme $T$ and a morphism $T \to \Branch^n(\mathbb{P}^n)$ corresponding to a $T$-family of branchvarieties $f: X \to \mathbb{P}^n_T$. The fiber product $\pr^{-1}(T):= T \times_{\Branch^n(\mathbb{P}^n)} \Branch^n(\bP^{N} \setminus \bP^{N - n - 1})$ is the functor that sends a $T$-scheme $T' \to T$ to the set of morphisms $\widetilde{f}: X_{T'} \to \mathbb{P}^N_{T'} \setminus \mathbb{P}^{N-n-1}_{T'} \cong \Tot_{\bP^n_{T'}}(\cO_{\bP^n_{T'}}(1)^{\oplus (N-n)})$ such that $\pr \circ \widetilde{f} = f_{T'}$. Such a morphism $\widetilde{f}$ amounts to equipping the $\cO_{\mathbb{P}_{T'}}$-algebra $f_*(\cO_X)_{T'}$ with the structure of a $\Sym^{\bullet}_{\cO_{\mathbb{P}_{T'}}}(\cO_{\mathbb{P}_{T'}}(-1)^{\oplus(N-n)})$-module. This, in turn, is equivalent to a morphism of $\cO_{\mathbb{P}^n_{T'}}$-modules $\cO_{\mathbb{P}^n_{T'}}(-1)^{\oplus (N-n)} \to f_*(\cO_X)_{T'}$. We conclude that $\pr^{-1}(T)$ is the functor $\Hom_{T}(\cO_{\mathbb{P}^n_{T}}(-1)^{\oplus (N-n)}, f_*(\cO_X))$ as in \cite[Lem. 2.16]{halpern2024moduli}, which is represented by a relatively affine scheme of finite type over $T$.

\end{proof}

\subsection{The \texorpdfstring{$\Chow$}{Chow} morphism} \label{subsection: chow morphism}

For clarity in this section, we introduce a vector space $U$ of dimension $N+1$, and regard $\bP^N = \bP(U) = \Proj(\Sym(U^\vee))$. We recall the definition of the Chow morphism for branchvarieties. For a given Hilbert polynomial $P(t) = \sum_{i=0}^n a_i \binom{t}{i}$, we set $d = a_{n-1} - n a_n$,
\[
V_{P,N} = \bigotimes_{i=1}^{n + 1} \Sym^d(U).
\]
Note that $V_{P,N} \cong \Gamma(\bP(U^\vee)^{n+1}, \cO(d,\ldots,d))$, so $\bP(V_{P,N})$ is identified with the scheme whose $T$-points are relative Cartier divisors in $\bP(U^\vee)^{n+1}_T$ of multidegree $(d,\ldots,d)$.

The Chow morphism is defined to be the composition
\[ \Chow: \Branch^P(\bP(U)) \to \Chow_{n, d}(\bP(U)) \to \bP(V_{P,N}) \]
The first morphism, defined in \cite[Remark 13.5]{rydh2008families}, sends a branchvariety $X \to \bP(U)_T$ to the pushforward of the fundamental cycle of $X$, which is equidimensional of dimension $n$ and of degree $d$. The second morphism sends a relative cycle to a corresponding relative Cartier divisor and, as noted in \cite[\S 17]{rydh2008families}, this relative Cartier divisor agrees with the one constructed by Mumford in \cite[\S 5.4]{mumford1994geometric}, and what is called the Chow divisor in \cite[Pg. 53]{knudsen-mumford}.

The paper \cite[Pg. 53]{knudsen-mumford} also describes the global section defining the Chow divisor more explicitly. 
For a Noetherian $k$-scheme $T$ and a morphism $f: T \to \Branch^P(\bP(U))$ corresponding to a branchvariety $X \to \bP(U)_T$, we set $b_{n+1}(X) := f^\ast b_{n+1}$. 
Let $\check{\pi}: \bP(U^\vee)^{n+1}_T \to T$ be the projection morphism, and note that 
\[ \check{\pi}_{\ast} \cO_{(\bP(U^\vee)^{n+1})_T}(d, \ldots, d) = V_{P,N} \otimes \cO_T. \]
The line bundle $b_{n+1}(X)$ is equal to the line bundle $\cM_{n + 1}$ discussed in \cite{knudsen-mumford}, so \cite[Thm. 4]{knudsen-mumford} shows that the Chow divisor is defined by a section $s(X) \in H^0(T, b_{n+1}(X) \otimes V_{P,N})$. Translating this to our setting gives the following:
\begin{prop} \label{prop: b_{n+1} is the pullback of O(1) under Chow}
We have $b_{n+1} = \Chow^\ast (\cO_{\bP(V_{P,N})}(1))$.
\end{prop}

Concretely, at any point $t\in T$, $s(X)$ gives a non-zero element of the fiber
\[
s(X)_t \in (V_{P,N})_{k(t)} \cong \Gamma(\bP(U^\vee)^{n+1}_{k(t)}, \cO(d,\ldots,d))
\]
that is well-defined up to a unit. For any collection of points $H_0,\ldots,H_n \in \bP(U^\vee)_{k(t)}$, the image of $X_t \to \bP(U)_t$ intersects $H_0 \cap \cdots \cap H_n \subset \bP(U)_t$ if and only if $s(X)_t(H_0,\ldots,H_n) = 0$; this is the defining property of the Chow form \cite[Defn. 8.9]{rydh2003chow}.

\begin{prop} \label{prop: Alg subfunctors are chow preimages of divisor complements}
For each linear subspace $L  \subseteq \bP(U)$ of codimension $n+1$,
\[ \Branch^P(\bP(U) \setminus L) = \Chow^{-1}(Y_L), \] where $Y_L \subset \bP(V_{P,N})$ is the complement of a hyperplane.
\end{prop}

\begin{proof}
As we have defined $\Branch(\bP(U) \setminus L)$ as an open substack of $\Branch(\bP(U))$, we may assume $k$ is algebraically closed,
in which case the claim amounts to an identification of sets of $k$-points. Suppose that $L$ is the common vanishing locus of $h_0, h_1, \ldots, h_n \in U^\vee$ and let $f: X \to \bP(U)$ be a branchvariety.

By definition $[X] \in \Branch(\bP(U) \setminus L)$ if $f(X) \cap L = \emptyset$, which as discussed above is equivalent to $s(X)(h_0,\ldots,h_n) \neq 0$. Under the identification of an element $s \in \Sym^d(U)$ with a function on $U^\vee$, one has $s(h) = (s,h^d)$, where $(-,-) : \Sym^d(U) \otimes \Sym^d(U^\vee) \to k$ is the canonical pairing $(u_1 \cdots u_d, h_1 \cdots h_d) = \frac{1}{d!} \sum_{\sigma \in S_d} h_1(u_{\sigma(1)}) \cdots h_d(u_{\sigma(d)})$.
Likewise, regarding $s(X) \in \Sym^d(U)^{\otimes(n+1)}$ as a function on $\bP(U^\vee)^{n+1}$ up to a unit, we have
\[
s(X)(h_0,\ldots,h_n) = (s(X),h_0^d\otimes \cdots \otimes h_n^d),
\]
where $(-,-)$ denotes the canonical pairing 
\[\Sym^d(U)^{\otimes(n+1)} \otimes \Sym^d(U^\vee)^{\otimes(n+1)} \to k.\]
 It follows that $[X] \in \Branch(\bP(U) \setminus L)$ if and only if $s(X)$ lies in the complement of the hyperplane in $\bP(V_{P,N})$ defined by $h_0^d \otimes \cdots \otimes h_n^d \in \Sym^d(U^\vee)^{\otimes(n+1)} \cong (\Sym^d(U)^{\otimes(n+1)})^\vee$.
\end{proof}

\begin{notn}
    Let $X$ be a projetive $k$-scheme equipped with a very ample line bundle $\cO_X(1)$ which induces an embedding $X \hookrightarrow \mathbb{P}^N$. For any given Hilbert polynomial $P$, we define the Chow morphism for $\Branch^P(X)$ to be the composition $\Chow: \Branch^P(X) \hookrightarrow \Branch^P(\mathbb{P}^N) \xrightarrow{\Chow} \mathbb{P}(V_{P,N})$.
\end{notn}

\section{(Semi)ample line bundles on \texorpdfstring{$\Branch$}{Branch}}

\subsection{Relative semiampleness of \texorpdfstring{$-b_n$}{-bn}} \label{subsection: semiample line bundles on branch}

\begin{prop} \label{prop: -b_n is relatively semiample}
Let $X$ be a projective $k$-scheme equipped with a very ample line bundle $\cO_X(1)$ which induces an embedding $X \hookrightarrow \mathbb{P}^N$. 
Fix a Hilbert polynomial $P$ of degree $n$.
Then the line bundle $-b_n$ is semiample relative to the Chow morphism $\Chow: \Branch^P(X) \to \bP(V_{P, N})$.
\end{prop}
\begin{proof}
The closed immersion $X \hookrightarrow \mathbb{P}^N$ induces a closed immersion of stacks $\Branch^P(X) \hookrightarrow \Branch^P(\mathbb{P}^N)$, and by definition the restriction of the line bundle $-b_n \in \Pic(\Branch^P(\mathbb{P}^N))$ to $\Branch^P(X)$ recovers the corresponding line bundle $-b_n \in \Pic(\Branch^P(X))$. Hence, it is sufficient to prove the result for $\Branch^P(\mathbb{P}^N)$, and so we may assume without loss of generality $X = \mathbb{P}^N$. We may restrict to the preimage $\Chow^{-1}(U)$ of a standard open $U \subset \bP(V_{P, N})$. By \Cref{prop: Alg subfunctors are chow preimages of divisor complements}, we have $\Branch^P(\bP^{N} \setminus \bP^{N - n - 1}) = \Chow^{-1}(U)$ for some linear subspace $\bP^{N - n - 1} \subset \bP^{N}$. Using the linear projection morphism from \Cref{lemma: linear projection}, we reduce to the case when the degree of the Hilbert polynomial $n$ matches the dimension of the target projective space $\mathbb{P}^n$. In this case, the semiampleness of $-b_n$ was shown in \Cref{coroll: discriminant morphism branch}.
\end{proof}

\subsection{Relative ampleness of \texorpdfstring{$-ab_n -\lambda_r$}{-abn-lambdar}}
\begin{lem} \label{lemma: normalization family branch constant}
    Suppose that the ground field $k$ is algebraically closed. Let $B$ be a smooth proper $k$-scheme, and choose $\varphi: B \to \Branch^n(\mathbb{P}^n)$ corresponding to a family of branchvarieties $f:X \to \mathbb{P}^n_B$. Suppose that $\varphi^*(2b_n)$ is the trivial line bundle on $B$. Then, after perhaps replacing $B$ with a finite \'etale cover, the normalization $X^{\nu}$ is a constant family of branchvarieties, i.e. there exists a normal scheme $Y$ equipped with a finite morphism $Y \to \mathbb{P}^n_{B}$ and an isomorphism of $\mathbb{P}^n_B$-schemes $X^{\nu} \cong Y \times B$.
\end{lem}
\begin{proof}
    By passing to the connected components of $B$, we may assume that $B$ is connected. \Cref{prop: discriminant morphism algred} and the assumption that $2b_n|_B$ is trivial imply that the associated discriminant divisor $D \subset \mathbb{P}^n_B$ of $f_*(\cO_X)$ (\Cref{defn: discriminant divisor}) is a constant divisor of the form $H \times B$ for some $H \subset \mathbb{P}^n$. Let $Z \subset \mathbb{P}^n_B$ denote complement of the maximal open subscheme where $f_\ast(\cO_X)$ is locally free. $Z$ has codimension at least $2$ in every fiber over $B$. 
    Recall that $f: X \to \mathbb{P}^n_B$ is \'etale on the open complement $W: = \mathbb{P}^n_B \setminus (D \cup Z)$ of the union $D \cup Z \subset \mathbb{P}^n_B$. Since $B$ is regular, the scheme $\mathbb{P}^n_B$ is regular. Zariski-Nagata purity \cite[X, Cor. 3.3]{SGA1} implies that $\pi_1(W) = \pi_1(\mathbb{P}^1_B \setminus D)$. On the other hand, since $D = H \times B$, we have $\pi_1(\mathbb{P}^1_B \setminus D) = \pi_1((\mathbb{P}^n \setminus H) \times B) = \pi_1(\mathbb{P}^n \setminus H) \times \pi_1(B)$, where the last equality follows because $B$ is proper and connected \cite[X, Cor. 1.7]{SGA1}. We conclude that $\pi_1(W) = \pi_1(\mathbb{P}^n \setminus H) \times \pi_1(B)$.
    
    The restriction $f: f^{-1}(W) \to W$ is a finite \'etale cover of $W$. The equality $\pi_1(W) = \pi_1(\mathbb{P}^n \setminus H) \times \pi_1(B)$ implies that, after perhaps replacing $B$ with a finite \'etale cover, there is a finite \'etale morphism $Y \to \bP^N \setminus H$ such that $f: f^{-1}(W) \to W$ is the restriction of $Y \times B \to (\mathbb{P}^n \setminus H) \times B$ over $W \subset(\mathbb{P}^n \setminus H) \times B$. Let $\widetilde{Y} \to \mathbb{P}^n$ denote the relative normalization of the composition $Y \to \mathbb{P}^n \setminus H \hookrightarrow \mathbb{P}^n$ \cite[\href{https://stacks.math.columbia.edu/tag/0BAK}{Tag 0BAK}]{stacks-project}. The normalization of the original family $X^{\nu}$ and the scheme $\widetilde{Y}\times B$ are both normal schemes. We have finite morphisms $X^{\nu} \to \mathbb{P}^n_B$ and $\widetilde{Y}\times B \to \mathbb{P}^n_B$, both of which are dominant when restricted to every irreducible component of the source. Furthermore, by construction, both schemes are isomorphic over the open $W \subset \mathbb{P}^n_B$. We conclude that $X^{\nu}$ and $\widetilde{Y} \times B$ are isomorphic $\mathbb{P}^n_B$-schemes (indeed, they are both forced to be isomorphic to the relative normalization of the composition $f^{-1}(W) \to W \hookrightarrow \mathbb{P}^n_B$).
\end{proof}

\begin{prop} \label{prop: -lambda_n relatively ample discriminant}
    Fix a Hilbert polynomial $P$ of degree $n$. Then, for all sufficiently large integers $r \gg 0$, the line bundle $-\lambda_r$ is ample relative to the discriminant morphism $\disc: \Branch^P(\mathbb{P}^n) \to \mathbb{P}\left(H^0(\mathbb{P}^n, \cO(d))\right)$ (as in \Cref{coroll: discriminant morphism branch}).
\end{prop}
\begin{proof}
    To show the relative ampleness of $-\lambda_r$, we may replace $\Branch^P(\mathbb{P}^n)$ with a closed fiber $\mathcal{Z} \subset \Branch^P(\mathbb{P}^n)$ of $\disc$. By the Nakai-Moishezon criterion \cite[Thm. 3.11]{kollar_projectivity} applied to the coarse space of $\cZ$ jointly with \cite[Thm. B]{rydh-approximation}, it suffices to show that for all proper integral algebraic spaces $B$ equipped with a generically finite morphism $\varphi: B \to \mathcal{Z}$, the top self-intersection of $\varphi^*(-\lambda_r)$ is positive for some uniform bound $r \gg 0$. 
    
    Choose such a morphism $\varphi: B \to \mathcal{Z}$. By Chow's lemma \cite[\href{https://stacks.math.columbia.edu/tag/088U}{Tag 088U}]{stacks-project} and resolution of singularities \cite{resolution_singularities}, we may assume that $B$ is a projective smooth $k$-scheme. Since $-2b_n$ is pulled back under the discriminant morphism $\disc$ (\Cref{coroll: discriminant morphism branch}), we see that $2b_n$ is the trivial line bundle on $\mathcal{Z}$. Hence, the morphism $\varphi$ corresponds to a $B$-family of branchvarieties $X \to \bP^n_B$ where $\varphi^*(2b_n)$ is trivial. By \Cref{lemma: normalization family branch constant}, after perhaps replacing $B$ with a finite \'etale cover, the normalization $X^{\nu}$ is a constant $B$-family of branchvarieties $Y \times B \to \mathbb{P}^n_B$. There is an injection $\iota: \cO_X \hookrightarrow \cO_{X^{\nu}}$, which remains injective after passing to $b$-fibers, since $X^{\nu}_b \to X_b$ is a dominant morphism of reduced schemes for all $b \in B$. Since $\cO_{X^{\nu}} = \cO_{Y \times B}$ is $B$-flat, the slicing criterion for flatness \cite[\href{https://stacks.math.columbia.edu/tag/046Y}{Tag 046Y}]{stacks-project} implies that the cokernel of $\cO_X \hookrightarrow \cO_{X^{\nu}}$ is $B$-flat. Therefore, we get a morphism $\psi: B \to \Quot^P(\cO_Y)$ to the quot scheme of subsheaves of the coherent sheaf $\cO_Y$ on $\mathbb{P}^n$ that have Hilbert polynomial $P$ and whose cokernel is flat over the base. Observe that the generic finiteness of $\varphi: B \to \mathcal{Z}$ implies that $\psi$ is generically finite: indeed, the normalization $X^{\nu}$ is determined from $X$ up to a finite group of isomorphisms and the algebra structure of $X$ can be recovered from that of $\cO_{X^{\nu}}$ via the subsheaf inclusion $\cO_X \hookrightarrow \cO_{X^{\nu}}$. Since the polynomial family of line bundles $-\lambda_r$ associated to the universal subsheaf is ample on the Quot scheme $\Quot^P(\cO_Y)$ for all sufficiently large $r\gg0$ \cite[Prop. 2.2.5]{huybrechts2010geometry}, it follows that the top-self intersection of $\psi^*(-\lambda_r)$ is positive. By construction $\psi^*(-\lambda_r) = \varphi^*(-\lambda_r)$, and therefore we conclude the desired positivity of the top self-intersection of $\varphi^*(-\lambda_r)$. 

To see that the bound on $r$ can be made uniform, we may stratify $\mathcal{Z} = \sqcup_i \mathcal{Z}_i$ by finitely many locally closed substacks where the normalization of $X$ is flat and its formation commutes with passing to geometric fibers. The sheaf $\cO_Y$ appearing in our argument must belong to the bounded family parametrized by $\sqcup_i \mathcal{Z}_i$, and therefore we may bound uniformly the necessary $r \gg 0$ that guarantees ampleness of $-\lambda_r$ for any Quot scheme $\Quot^P(\cO_Y)$ as in the proof.
\end{proof}

\begin{prop} \label{prop: -lambda_n is relatively ample}
Let $X$ be a projective $k$-scheme equipped with a very ample line bundle $\cO_X(1)$ which induces an embedding $X \hookrightarrow \mathbb{P}^N$. 
Fix a Hilbert polynomial $P$ of degree $n$.
Then for all sufficiently large integers $r \gg 0$, there exists $a(r)>0$ such that for all rational numbers  $a \geq a(r)$ the $\mathbb{Q}$-line bundle $-ab_n -\lambda_r$ is ample relative to the Chow morphism $\Chow: \Branch^P(X) \to \bP(V_{P, N})$.
\end{prop}
\begin{proof}
We immediately reduce to the case when $X = \mathbb{P}^N$ and the ground field $k$ is algebraically closed. We may work over the standard affine open subsets $U \subset \bP(V_{P, N})$. By \Cref{prop: Alg subfunctors are chow preimages of divisor complements}, we have $\Branch^P(\bP^{N} \setminus \bP^{N - n - 1}) = \Chow^{-1}(U)$ for some linear subspace $\bP^{N - n - 1} \subset \bP^{N}$. Since the ampleness of a line bundle is preserved by pullback under affine finite type morphisms, the affine linear projection morphism from \Cref{lemma: linear projection} allows us to reduce to the case when the degree of the Hilbert polynomial $n$ equals the dimension of the target projective space $\mathbb{P}^n$. In this case we have that the target of the Chow morphism is $\Spec(k)$, and $\Branch^P(\mathbb{P}^n) \cong \Alg^P(\mathbb{P}^n)_{\reduced}$ (\Cref{lemma: branchn isomorphic to algred}). Since the line bundle $-2b_n$ is the pullback of an ample line bundle under the discriminant morphism $\disc: \Branch^P(\mathbb{P}^n) \to \mathbb{P}\left(H^0(\mathbb{P}^n, \cO(d))\right)$ from \Cref{coroll: discriminant morphism branch}, to conclude it suffices to show that $-\lambda_r$ is ample relative to $\disc$ for sufficiently large $r\gg0$. This is the content of \Cref{prop: -lambda_n relatively ample discriminant}. 
\end{proof}
\subsection{Nefness of \texorpdfstring{$nb_{n+1}-b_n$}{nbn+1 - bn}}

\begin{lem} \label{lemma: slicing branchvarieties}
    Suppose that $k$ is algebraically closed. Let $C$ be a smooth projective connected curve over $k$, and let $f:X \to \bP^N_C$ be a $C$-family of branchvarieties in $\Branch^n(\bP^N_k)(C)$. If $N>n$, then there exists a codimension $n$ linear subspace $\bP(V) \subset \bP^N$ such that the fiber product $Y:= f^{-1}(\bP(V)_C)$ satisfies that $Y \to C$ is finite, flat and with reduced generic fiber.
\end{lem}
\begin{proof}
    By a dimension counting argument, there exists a linear subspace $\bP(W) \subset \bP^N$ of dimension $n+1$ such that $f(X)$ is contained in the locus of definition of the rational linear projection $\bP^N_C \dashrightarrow \bP(W)_C$. We can view the composition with the linear projection $X \to \bP(W)_C$ as a $C$-family of branchvarieties in $\Branch^n(\bP(W))(C)$, and so we may assume without loss of generality that $N=n+1$ for the rest of the proof.
    
    The scheme $X$ is equidimensional and reduced, since the flat morphism $X \to C$ has reduced equidimensional fibers. By induction on $n$, it suffices to show the following:

    \noindent \textbf{Claim:} For $n \geq 1$, let $f: X \to \bP^{n+1}_C$ be a finite morphism such that 
    \begin{enumerate}
        \item[(A)] $X$ is reduced and equidimensional of dimension $n+1$.
        \item[(B)] The projection $g:X \to C$ is flat.
    \end{enumerate}
    Then there exists a hyperplane $H \subset \bP^{n+1}$ such that the preimage $Y:=f^{-1}(H_C)$ satisfies the following:
    \begin{enumerate}
        \item[(A')] $Y$ is reduced and equidimensional of dimension $n$.
        \item[(B')] $Y \to C$ is flat.
    \end{enumerate}

    We conclude the proof by proving the claim. By a general form of Bertini's theorem (see, for example, \cite[Thm. 1.4]{bertini-theorems-tanaka}), for a general hyperplane $H \subset \bP^{n+1}$ we have that $Y = f^{-1}(H_C)$ satisfies (A'). For any such general hyperplane $H$, the morphism $Y \to C$ to the smooth curve $C$ is flat if and only if $Y \to C$ has fibers of dimension $n-1$. Since the fibers of $g: X \to C$ have dimension $n$, this holds if and only if we have the following:
    
    \noindent (*) For all points $c \in C(k)$ the preimage $f^{-1}(H_C)$ does not contain any irreducible component of $g^{-1}(c) \subset X$.

    We show that a general hyperplane $H$ satisfies (*) by a dimension counting argument. First, note that the scheme theoretic image $Z:= f(X) \subset \bP^{n+1}_C$ is a divisor such that the fibers of $g: Z \to C$ are equidimensional of dimension $n$, and therefore it is a $C$-flat relative effective Cartier divisor on $\bP^{n+1}_C \to C$. We may replace $X \to \bP^{n+1}_C$ with $Z \hookrightarrow \bP^{n+1}_C$ in order to show (*). The $C$-family of divisors $Z \hookrightarrow \bP^{n+1}_C$ corresponds to a morphism $h: C \to \bP(H^0(\bP^{n+1}_k, \cO(m)))$ for some $m \geq 1$. To simplify notation, let us denote $S_m:= H^0(\bP^{n+1}_k, \cO(m))$. We have a morphism $\mu: \bP(S_1) \times \bP(S_{m-1}) \to \bP(S_m)$ induced by multiplication of homogeneous polynomials. Consider the fiber product 
   \[\begin{tikzcd}
				F \ar[r] \ar[d] & \bP(S_1) \times \bP(S_{m-1}) \ar[d, "\mu"]\\   C \ar[r, "h"] & \bP(S_m)
			\end{tikzcd}\]
   The general hyperplane $H \subset \bP^{n+1}$ satisfies (*) if and only if the projection $F \to \bP(S_1) \times \bP(S_m) \to \bP(S_1)$ is not dominant. Note that the fibers of $F \to C$ are finite, by unique factorization of polynomials. Therefore $\dim(F) \leq 1$. On the other hand, we have $\dim(\bP(S_1)) = n+1 \geq 2$. Hence, $F \to \bP(S_1)$ cannot be dominant, as desired.
\end{proof}

\begin{lem} \label{lemma: slicing branchvarieties when n=N}
    Suppose that $k$ is algebraically closed. Let $C$ be a smooth projective connected curve over $k$, and let $f:X \to \bP^N_C$ be a $C$-family of branchvarieties in $\Branch^n(\bP^N_k)(C)$. If $N=n$, then there exists a point $p \in \bP^N(k)$ such that the fiber product $Y:= f^{-1}(p \times C)$ satisfies that $Y \to C$ is finite, flat and with reduced generic fiber.
\end{lem}
\begin{proof}
    The pushforward $f_*(\cO_X)$ is a torsion-free sheaf on $\bP^N_C$. In particular, it is a vector bundle away from a closed subset $Z \subset \bP^N_C$ of codimension at least $2$. The composition $Z \to \bP^N_C \to \bP^N$ is not dominant for dimension reasons, and therefore for a general point $p \in \bP^N(k)$ we have that $p \times C$ is disjoint from $Z$. It follows that for all such general $p$ the fiber product $Y: = f^{-1}(p \times C)$ is finite and flat over $C$, and in particular it is $S_1$ and equidimensional of dimension $1$. For any given fixed point $c \in C(k)$, a generic point $p$ will lie in the nonempty open locus of $\bP^N$ over which the morphism of reduced schemes $f:X_c \to \bP^N \times c \cong \bP^N$ is \'etale. For such a generic point, the flat morphism $Y = f^{-1}(p \times C) \to C$ will be unramified at $c$, and hence $Y \to C$ will be \'etale over $c$. It follows that there is an open subset $U \subset C$ such that $Y$ is reduced over the open dense preimage $Y_U$. Since $Y$ is $S_1$, it follows that $Y$ is reduced, and this in turn implies that the generic fiber of $Y \to C$ is reduced.
\end{proof}

\begin{prop} \label{prop: -bd nef on branch}
    Let $X$ be a projective $k$-scheme equipped with a very ample line bundle $\cO_X(1)$. The line bundle $n b_{n+1} -b_n$ on $\Branch^n(X)$ is nef.
\end{prop}
\begin{proof}
    We may immediately reduce to the case when $X= \bP^N$ and $k$ is algebraically closed. Let $C$ be a smooth projective curve over $k$, and choose a morphism $\varphi: C \to \Branch^n(\bP^N)$ corresponding a $C$-family $f:Z \to \bP^N_C$ of equidimensional branchvarieties. By \Cref{lemma: slicing branchvarieties} or \Cref{lemma: slicing branchvarieties when n=N}, there exists some linear subspace $\bP^{N-n} \subset \bP^N$ such that the induced morphism $g: Y:= f^{-1}(\bP^{N-n}_C) \to C$ is finite and flat. Using the formula for $b_n \in \Pic(\Branch^n(\bP^N))$ in \Cref{remark: taylor series expansion}, it can be checked that we have an isomorphism $\varphi^*(nb_{n+1} - b_n) \cong \det(g_*(\cO_Y))^{-1}$.\footnote{
    For simplicity, given a morphism $\psi:C \to \Branch^n(\bP^N)$ corresponding to a $C$-family $W \to \bP^N_C$ we let $\lambda_r(W) := \varphi^{\ast}(\lambda_r)$. In particular, we need to show that $nb_{n+1} - b_n \simeq -\lambda_0(Y)$. 
    Choose hyperplanes $H_1,\dots, H_n \subseteq \bP^N_C$ such that $\bP^{N-n}_C = H_1 \cap \dots \cap H_n$. If we let $X_i = f^{-1}(H_1 \cap \cdots \cap H_i)$, then for any $r \in \mathbb{Z}$ we have exact sequences 
    $$
    0 \to \cO_{X_i}(r-1) \to \cO_{X_i}(r) \to \cO_{X_{i+1}}(r) \to 0.
    $$
    From this we obtain $\lambda_r(X_{i+1}) \simeq \lambda_r(X_i) - \lambda_{r-1}(X_i)$. In particular, $$\lambda_r(X_{i+1}) \simeq \sum_{j=0}^{i+1}(-1)^j\binom{i+1}{j}\lambda_{r-j}(X).$$
    Choosing $i+1 = n$ and $r=0$, we can use the formula in \Cref{remark: taylor series expansion} and rearrange to obtain
    \begin{align*}
    \lambda_0(Y) = \lambda_0(X_n) \simeq \sum_{j=0}^n (-1)^j \binom{n}{j}\lambda_{-j}(X)
    = \sum_{\ell=0}^{n+1} \left(\sum_{j=0}^n(-1)^{\ell+j}\binom{n}{j}\binom{j+\ell-1}{\ell}\right)b_{\ell}.
    \end{align*} 
    Most of these inner sums vanish and we obtain $\lambda_0(X_n) \simeq b_n-nb_{n+1}$. In particular, $-\lambda_0(X_n) \simeq nb_{n+1}-b_n$, as required.
    } It suffices to show that $\deg(\det(g_*(\cO_Y))^{-1}) \geq 0$. 
    
    Let $s: \cO_{C} \to \det(g_*(\cO_Y))^{-2}$ denote the discriminant section of the flat $\cO_C$-algebra $g_*(\cO_Y)$ (\Cref{notn: discriminant}). Since the generic fiber of $g:Y \to C$ is reduced and the characteristic of $k$ is $0$, it follows that $Y \to C$ is \'etale over an open subscheme of the target $C$. It follows from \Cref{notn: discriminant} that $s : \cO_C \to \det(g_*(\cO_Y))^{-2}$ is nonzero, and therefore we must have $\deg(\det(g_*(\cO_Y))^{-2}) \geq 0$. We conclude that $\deg(\det(g_*(\cO_Y))^{-1}) \geq 0$, as desired. 
    
\end{proof}

\subsection{Semiampleness of \texorpdfstring{$ab_{n+1} - b_n$}{abn+1 - bn}}

\begin{prop} \label{prop: semiampleness of a b_n+1 - b_n}
      Fix a Hilbert polynomial $P$ of degree $n$. Let $X$ be a projective $k$-scheme equipped with a very ample line bundle $\cO_X(1)$. Then for every positive rational number $a>n$ the $\mathbb{Q}$-line bundle $ab_{n+1} - b_n$ is semiample on $\Branch^P(X)$.
\end{prop}
\begin{proof} 
We may reduce immediately to the case when $X= \bP^N$ and $k$ is algebraically closed. Let $ \Branch^P(\bP^N) \to M$ denote the projective coarse space. After replacing $b_{n+1}$ and $b_n$ by a common multiple, we may assume that $b_{n+1}$ and $b_n$ descend to $M$, and we are reduced to checking semiampleness on $M$. Since $-b_n$ is relatively semiample for $\Chow: \Branch^P(\bP^N)\to \bP(V_{P,N})$ (\Cref{prop: -b_n is relatively semiample}), we may further assume that there is a projective scheme $g: Y \to \bP(V_{P,N})$, a $\bP(V_{P,N})$-ample line bundle $H$ on $Y$ and a morphism $f: M \to Y$ such that $-b_n \cong f^*(H)$. Recall that $b_{n+1}$ is the pullback
of a very ample line bundle $\cO_{\bP(V_{P,N})}(1)$ on the Chow variety $\bP(V_{P,N})$ by \Cref{prop: b_{n+1} is the pullback of O(1) under Chow}. By abuse of notation, we also denote $\cO_{\bP(V_{P,N})}(1)$ by $b_{n+1}$. There is a sufficiently large positive integer $m\gg0$ such that $m \cdot g^*(b_{n+1}) +H$ is ample on $Y$. We may assume that $m>a$. Since $f^*(ng^*(b_{n+1}) + H) = n b_{n+1} -b_n$ is nef on the coarse space $M$ (\Cref{prop: -bd nef on branch}), it follows that $ng^*(b_{n+1}) + H$ is nef on $Y$. It follows that for every nonnegative rational number $l \geq 0$, the rational line bundle 
\[(m  + ln) \cdot g^*(b_{n+1}) + (1+l)H = (m \cdot g^*(b_{n+1}) + H) +l(n g^\ast(b_{n+1}) + H)\]
is ample on $Y$, because it is the sum of an ample and a nef rational line bundle. Hence, $f^*((m  + ln +1) \cdot g^*(b_{n+1}) +(1+l) H) = (m  + ln +1)b_{n+1} - (1+l) b_n$ is semiample for all $l \geq 0$. Setting $l = \frac{m-a}{a-n}$ we get that $(m+ln)b_{n+1} - (1+l) b_n$ is a positive multiple of the rational line bundle $ab_{n+1} - b_n$, and we therefore conclude that $ab_{n+1} - b_n$ is semiample.
\end{proof}

\subsection{Some other relevant line bundles}

In \cite[Rem.~5.2]{branchvarieties-paper}, the question is raised as to whether the line bundle $\lambda_{m\ell} - (n-1/2)m \lambda_{\ell}$ is ample on $\Branch$ when $m,\ell \gg 0$. This turns out to be false, because the leading term is not nef. For instance, if one considers $\Branch^n(\bP^n)$ and $n>1$, then
\[
\lambda_{m\ell} - (n-1/2)m \lambda_{\ell} = \frac{(m\ell)^n-(n-1/2)m\ell^n}{n!} b_n + O((m\ell)^{n-1})
\]
If one chooses a curve in $\Branch^n(\bP^n)$ on which the semiample line bundle $-b_n$ has positive degree, such as one whose generic point corresponds to a variety satisfying condition $S_2$ (see \Cref{thm: ampleness of cm line bundle on s_2 locus} below), then the leading term above will have negative degree for all $m\geq n$ and $\ell>0$.

\medskip

Another natural sequence of line bundles associated to a family of polarized varieties is $m P(m) \lambda_{\ell} - \ell P(\ell) \lambda_m$ for $0 \ll m \ll \ell$.
The weight of this sequence is used to define various intrinsic notions of stability for polarized varieties \cite{ross-thomas}. Using the expansion of $\lambda_\ell$ given in \Cref{defn: knudsen-mumford coefficients}, one can compute the leading order term of this sequence, c.f., \cite[Prop. 3.1]{fine2006note}. The leading order term is a positive multiple of the following $\bQ$-line bundle:

\begin{defn} \label{defn: CM line bundle}
The CM-line bundle  $\CM \in \Pic(\Branch^n(\bP^N))_\bQ$ is defined to be
$$
\CM := \left(n(n+1) + 2\frac{a_{n-1}'}{a_{n}'}\right) b_{n+1} -2(n+1)b_n.
$$
Here we denote by $a_n'$ and $a_{n-1}'$ the top coefficients in the monomial expansion $P(t) = \sum_{i=0}^n a_i' t^i$ of the Hilbert polynomial of the branchvariety.
\end{defn}

One challenge in constructing moduli spaces for polarized varieties is that this line bundle is not semiample on Hilbert schemes \cite{fine2006note}, and therefore not suitable for applying GIT to Hilbert schemes. We observe that this remains the case for $\Branch^P(\mathbb{P}^  N)$ in general.

\begin{ex}
    Let $T = \bP^1$ and consider the vector bundle $\mathcal{E} = \cO_{T}(2) \oplus \cO_{T}(-1)^{\oplus 2}$ on $T$. 
    Define the threefold $X' = \bP(\mathcal{E}) = \Proj(\Sym(\mathcal{E}^\vee))$ and consider the curve $C = \bP(\cO_{T}(2)) \subseteq X'$ corresponding to the inclusion of the summand $\cO(2) \hookrightarrow \mathcal{E}$. Let $X = \mathrm{Bl}_CX'$ and let $\pi$ denote the composition of the blowup $q:X \to X'$ and the natural projection $X' \to T$. Then $\pi:X \to T$ is a smooth family, where for each $t\in T$, the fiber $X_t$ is isomorphic to a blowup of $\bP^2$ at a point, a smooth del Pezzo surface of degree $8$. In fact, the family $X \to T$ is trivializable over the complement of any point in $T$. If $E$ denotes the exceptional divisor of the blowup $q:X \to X'$, then the $\bQ$-line bundle $\cL_{\varepsilon} := q^{\ast}(\cO_{X'}(1))\otimes \cO(-\varepsilon E)$ for $0 < \varepsilon \ll 1$ is $\pi$-relatively ample. It was shown in \cite[Ex. 5.3]{fine2006note} that $\deg(\lambda_{\text{CM}}(\cL_\varepsilon)) < 0$, and explored further in \cite[Ex. 12.1]{codogni2021positivity}. To obtain a polarization (as opposed to a $\bQ$-polarization) for the family $X \to T$, we fix an $a>0$ sufficiently divisible so that $\cL^a_\varepsilon$ is a line bundle, and $\deg(\lambda_{\text{CM}}(\cL_\epsilon^a)) = a^2 \deg(\lambda_{\text{CM}}(\cL_\epsilon)) < 0$ by \cite[Rem.~2.6]{fine2006note}.

    We realize $X \to T$ as a family of branchvarieties by choosing a $b \gg 0$ such that $\pi_\ast(\cL_\varepsilon^a) \otimes \cO_{\bP^1}(b)$ is a globally generated vector bundle on $\bP^1$. Then a choice of surjection $\cO_T^N \twoheadrightarrow \pi_\ast(\cL_\varepsilon^a) \otimes \cO_{\bP^1}(b)$ for some $N \gg 0$ gives a closed immersion of $T$-schemes $f:X \to \bP^N_T$ such that $f^\ast(\cO(1)) \cong \cL_{\varepsilon}^a \otimes \pi^{\ast}(\cO_{\bP^1}(b))$. In particular, $X$ is now a family of branchvarieties in $\bP^N$, corresponding to a morphism $\varphi: T \to \Branch^P(\bP^N)$. By \cite[Cor. 4.6]{fine2006note}, tensoring by a line bundle from the base does not affect $\lambda_{\text{CM}}$, so
    \begin{align*}
    \varphi^{\ast}(\CM) 
        = \lambda_{\text{CM}}(X,\cL_\varepsilon^a \otimes \pi^{\ast}(\cO_{\bP^1}(b)))
        = \lambda_{\text{CM}}(X,\cL_\varepsilon^a) 
       = a^2\lambda_{\text{CM}}(X,\cL_\varepsilon).
    \end{align*}
has negative degree on $T$. This shows that $\CM$ is not nef on $\Branch$.
\end{ex}

The counterexample above is very robust. It suggests that, although $\cL_{\text{CM}}$ is known to be ample on the moduli space of $K$-semistable klt Fano varieties \cite{codogni2021positivity, chenyang_positive}, this result can not be obtained from applying GIT to any variant of the Hilbert scheme. More precisely, $\CM$ does not extend to a nef line bundle on \emph{any} projective variety that compactifies the quasi-projective variety of smooth surfaces embedded in $\bP^N$. Although one can use $\CM$ to define a notion of $K$-semistability for arbitrary polarized varieties, it might be that more sophisticated notions of semistability are needed to construct moduli spaces outside of the general type and Fano cases.

In some sense, though, $\CM$ is closer to being ample on $\Branch^P(\bP^N)$ than on the Hilbert scheme. For instance, $\CM$ is semiample on $\Branch^P(\bP^N)$ in the following special cases:
\begin{itemize}
    \item If $a_{n-1}'/a_n' > n (n+1)/2$, then $\CM$ is semiample on $\Branch^P(\bP^N)$ by \Cref{prop: semiampleness of a b_n+1 - b_n}.
    \item If $N=\deg(P)=n$, then $\CM$ is semiample, because in this case the line bundle $b_{n+1}$ is trivial.
\end{itemize}

\section{Positivity of \texorpdfstring{$-b_n$}{-bn} on the deminormal locus}
This section is dedicated to proving \Cref{thm: ampleness of cm line bundle on s_2 locus}, which is an ampleness statement when we restrict to families with deminormal singularities.

\begin{defn}[{\cite[Defn. 11.10]{kollar-families-varieties-book}}]
    Let $X$ be a scheme of finite type over a field $K \supset k$. We say that $X$ is \textbf{deminormal} if it is equidimensional, it satisfies Serre's condition $S_2$ and it has at worst nodal singularities in codimension $1$.
\end{defn}

\begin{defn}
    Given any fixed Hilbert polynomial $P$, we denote by $\Branch^P_{\mathrm{Dn}}(X)$ the substack of $\Branch^P(X)$ parametrizing families of branchvarieties $f: Z \to X \times T$ such that every fiber of $Z \to T$ is deminormal. This is an open substack by \cite[Cor. 10.42]{kollar-families-varieties-book}.
\end{defn}

We will need the following technical lemma.
\begin{lem} \label{lemma: constant family deminormal}
    Suppose that $k$ is algebraically closed. Let $C$ be a smooth connected quasi-projective curve over $k$, and let $C \to \Branch^n(\mathbb{P}^n)$ be a morphism corresponding to a $C$-family of branchvarieties $X \to \mathbb{P}^n_C$. Suppose that the following are satisfied:
    \begin{enumerate}
        \item The fiber $X_{\eta}$ over the generic point $\eta \in C$ is deminormal.
        \item The normalization $X^{\nu}$ equipped with its natural morphism $X^{\nu} \to \mathbb{P}^n_C$ is isomorphic to a constant family $Y \times C \to \mathbb{P}^n_C$, where $Y \to \mathbb{P}^n$ is a branchvariety over $k$.
        \item The discriminant divisor $D \subset \mathbb{P}^n_C$ (as in \Cref{defn: discriminant divisor}) is of the form $F \times C$ for some Cartier divisor $F \subset \mathbb{P}^n$.
    \end{enumerate}
    Then the $C$-family of branchvarieties $X \to \mathbb{P}^n_C$ is isomorphic to a constant family of the form $\widetilde{Y} \times C \to \mathbb{P}^n_C$ for some branchvariety $\widetilde{Y} \to \mathbb{P}^n$ over $k$.
\end{lem}
\begin{proof}
    Since the stack $\Branch^n(\mathbb{P}^n)$ is separated, in order to show that the family $X \to \mathbb{P}^n_C$ is constant we may replace $C$ with a nonempty open subscheme. After shrinking $C$, we may assume that $X$ is a deminormal scheme.

     Let $F_Y \subset Y$ denote the preimage of the Cartier divisor $F \subset \bP^{n}$ under $Y \to \bP^{n}$. We denote by $\bigcup_{l} F_l = (F_Y)_{\red}$ the decomposition of the reduced subscheme $(F_Y)_{\red}$ into irreducible Weil divisors $F_l \subset Y$.
     Note that there is an open $U \subset \bP^{n}$ whose complement has codimension at least $2$ and such that the restrictions $(F_l)_{Y_U}$ of the Weil divisors $F_l$ to $Y_U \subset Y$ are disjoint and smooth.
     
     We will momentarily work over the open $U \subset \mathbb{P}^n$. For ease of notation, set $X' = X_{U \times C}$, and $Y' = Y_U$. Similarly let $F' := F_U$, and let $F'_l := (F_l)_{Y_U}$. The variety $X'$ and its normalization $(X')^{\nu} = Y' \times C$ fit into a commutative diagram
     \begin{equation} \label{eqn: conductor push-out diagram}
     \begin{tikzcd}
        \Cond \ar[d] \ar[r, symbol = \hookrightarrow] & Y' \times C \ar[d] \\
        G \ar[r, symbol = \hookrightarrow] & X',
    \end{tikzcd}         
    \end{equation}
    where $\Cond$ is the conductor subscheme of $Y' \times C \to X'$, defined locally over $X'$ by the ideal $\cI_{\Cond} = \{z \in \cO_{Y' \times C} \, | \, z \cdot \cO_{Y' \times C} \subset \cO_{X'}\}$, and $\Cond \to G$ is a finite morphism to the scheme theoretic image $G \subset X'$ of $\Cond \to X'$. As explained in \cite{mo186650}, the diagram \eqref{eqn: conductor push-out diagram} is a pushout diagram of schemes. Let us reproduce the argument from \cite{mo186650} here for completeness. Note that one may reduce to the case when $X'$ is affine, and then the cocartesianness of the diagram \eqref{eqn: conductor push-out diagram} amounts to checking that the naturally induced morphism from the ring $\cO_{X'}$ to the pullback $A$ of the diagram of rings 
    \[\cO_{(X')^{\nu}} \leftarrow \cO_{\Cond} \to \cO_G\]
    is an isomorphism. The composition $\cO_{X'} \to A \to \cO_{(X')^{\nu}}$ is injective, and hence $\cO_{X'} \to A$ is injective. To show surjectivity, choose an element $(s,\overline{r}) \in A$ corresponding to a pair of elements $s \in \cO_{(X')^{\nu}}$ and $\overline{r} \in \cO_{G}$ that have the same image in $\cO_{\Cond}$. Fix a preimage $r \in \cO_{X'}$ of $\overline{r}$ under the quotient morphism $\cO_{X'} \twoheadrightarrow \cO_{G'}$. The difference $s-r$ inside the normalization $\cO_{(X')^{\nu}}$ is sent to $0$ under the quotient morphism $\cO_{(X')^{\nu}} \twoheadrightarrow \cO_{\Cond} = \cO_{(X')^{\nu}}/\cI_{\Cond}$, and hence $r-s$ is contained in the conductor ideal $\cI_{\Cond}$. Note that $\cI_{\Cond} \subset \cO_{X'}$ by its definition, and so it follows that $s-r \in \cO_{X'}$, which in turn implies that $s \in cO_{X'}$. By construction, the image of $s \in \cO_{X'}$ in $A$ recovers $(s, \overline{r})$, thus concluding the proof of surjectivity and the fact that \eqref{eqn: conductor push-out diagram} is cocartesian.
    
    Since $X'$ is deminormal, $\Cond \subset Y' \times C$ is a reduced divisor (see the discussion in \cite[5.2, pg. 189]{kollar-singularities-mmp}). Note that $\Cond \to X'$ maps to the locus where $X'$ is singular. Hence $\Cond$ is contained in the reduced preimage $(F'_{Y})_{\red} \times C \subset Y' \times C$ of the closed locus $F' \times C \subset U\times C$ over which $X' \to U \times C$ is not \'etale (note that $X$ is smooth at every point where $X' \to U \times C$ is \'etale). In view of the reducedness of $\Cond$, this forces $\Cond \subset Y' \times C$ to be a union $\sqcup_{i \in I} F_i' \times C$ of some of the irreducible components of the smooth Cartier divisor $(F'_Y)_{\red} \times C \subset Y' \times C$. In particular, $\Cond$ is a normal scheme.
    
    Let us set $R := \sqcup_{i \in I} F_i'$, so that $\Cond = R \times C$. Since the characteristic of $k$ is $0$, we may assume after further shrinking $U \subset \mathbb{P}^n$ (keeping the property that its complement has codimension at least $2$) that $R \to F'_{\red}$ is \'etale. By the discussion in \cite[5.2, pg. 189]{kollar-singularities-mmp}, the scheme $G$ is obtained by quotienting $\Cond = R \times C$ via an involution $\tau$. In this case $\tau$ is an automorphism of $R \times C$ as a $F'_{\red} \times C$-scheme, since it is compatible with the projections to $\mathbb{P}^n_C$. We may think of $\tau$ as a $F'_{\red}$-morphism $h: C\times F'_{\red} \to \underline{\Aut}(R/F'_{\red})$, where $\underline{\Aut}(R/F'_{\red}) \to F'_{\red}$ denotes the finite \'etale $F'_{\red}$-group scheme parametrizing automorphisms of the finite \'etale $F'_{\red}$-scheme $R$. Since $C$ over $k$ is (geometrically) connected, the morphism $h: C\times F'_{\red} \to \underline{\Aut}(R/F'_{\red})$ factors through a section $F'_{\red} \to \underline{\Aut}(R/F'_{\red})$. In other words, there is an involution $\overline{\tau}: R \to R$ over $F'_{\red}$ such that $\tau = \overline{\tau} \times id_C: R \times C \to R \times C$. Therefore, the quotient $G=\Cond/\tau$ is of the form $G = (R/\overline{\tau}) \times C$, and the corresponding quotient morphism is a product of the form $\Cond = R \times C \xrightarrow{t \times id_C} (R/\overline{\tau}) \times C = G$.
    
    Let $Q$ denote the pushout of $R/\overline{\tau} \xleftarrow{t} R \hookrightarrow Y'$. Since $\overline{\tau}: R \to R$ is a morphism of $U$-schemes by construction, the composition $R \to Y' \to U$ factors through a morphism $R/\overline{\tau} \to U$. Therefore, the map $Y' \to U$ induces a finite morphism $g:Q \to U$, where we recall that $U \subset \bP^N$ has complement of codimension at least $2$. The reflexive hull of $g_*\cO_Q$ extends uniquely to a coherent $S_2$-algebra over $\bP^N$ \cite[\href{https://stacks.math.columbia.edu/tag/0EBJ}{Tag 0EBJ}]{stacks-project}, whose relative spectrum we denote by $\widetilde{Q} \to \bP^N$. By construction, there is an isomorphism of $U \times C$-schemes $\widetilde{Q}_U \times C \cong X' = X_{U \times C}$.
    \footnote{Indeed, we have shown that the normalization pushout square for $X'$ in \Cref{eqn: conductor push-out diagram} agrees with the defining pushout square for $Q = \widetilde{Q}_U$ base-changed to $C$.}
    Since $X$ is $S_2$, it follows that $\cO_{X}$ is a reflexive sheaf on $\bP^N \times C$ by \cite[\href{https://stacks.math.columbia.edu/tag/0AY6}{Tag 0AY6}]{stacks-project}. Hartogs' theorem \cite[\href{https://stacks.math.columbia.edu/tag/0EBJ}{Tag 0EBJ}]{stacks-project} then implies that $X \cong \widetilde{Q} \times C$. We conclude that $X$ is isomorphic to the base change to $C$ of the branchvariety $\widetilde{Q} \to \bP^N$, as desired.
\end{proof}

\begin{lem} \label{lemma: relative positivity of -b_n on branch-chow over s_2 locus}
    Let $C$ be a smooth projective connected curve over $k$, and let $\varphi: C \to \Branch^P(\bP^N)$ be a morphism satisfying the following conditions.
    \begin{enumerate}
        \item The set-theoretic image $\varphi(C) \subset |\Branch^P(\bP^N)|$ is not a single point.
        \item Image of the generic point of $C$ under $\varphi$ lands in $\Branch^P_{\mathrm{Dn}}(\bP^N)$.
        \smallskip
        \item The image $\varphi(C)$ is contained in a fiber of $\Chow: \Branch^P(\bP^N) \to \bP(V_{P,N})$.
    \end{enumerate}
    Then the pullback $\varphi^*(b_n^{-1})$ is ample on $C$.
\end{lem}
\begin{proof}
   We immediately reduce to the case when $k$ is algebraically closed. The family $\varphi$ corresponds to a flat family of reduced equidimensional schemes $\pi: X \to C$ equipped with a finite morphism $f: X \to \bP_C^N$. Since the image $\varphi(C)$ is entirely contained in a fiber of the Chow morphism, it is contained in the preimage $\Chow^{-1}(U)$ of one of the standard affine open subschemes as in \Cref{prop: Alg subfunctors are chow preimages of divisor complements}. Because the linear projection morphism of \Cref{lemma: linear projection} is affine, the composition $\widetilde{\varphi}: C \to \Chow^{-1}(U) \to \Branch^P(\bP^n)$ is also non-trivial, so it suffices to prove the lemma when $N=n$ is the degree of $P$. 
   Let $D \subset \mathbb{P}^n_C$ denote the discriminant divisor of $f_*(\cO_X)$ as in \Cref{defn: discriminant divisor}. By \Cref{prop: discriminant morphism algred}, $D$ is a relative effective Cartier divisor cut out by a section of $\cO(m) \boxtimes \varphi^*(b_n^{-2})$ on $\mathbb{P}^n \times C$ for some $m \geq 0$. Hence, in order to show that $\varphi^*(b_n^{-1})$ is ample, it is enough to prove that $D$ is not a horizontal divisor on $\bP_C^{n}$, i.e., $D$ is not of the form $F \times C$ for some Cartier divisor $F \subset \bP^n$.

But indeed, suppose for the sake of contradiction that $D = F \times C$ for some Cartier divisor $F \subset \bP^{n}$. By \Cref{lemma: normalization family branch constant}, after perhaps replacing $C$ with a finite \'etale cover, the normalization $X^{\nu} \to \mathbb{P}^n_C$ of $X$ is isomorphic to a constant family $Y \times C \to \mathbb{P}^n_C$ of branchvarieties, where $Y \to \mathbb{P}^n$ is a normal branchvariety over $k$. In view of this, by \Cref{lemma: constant family deminormal} it follows that the $C$ family $X \to \mathbb{P}^n_C$ is isomorphic to a constant family, contradicting our assumption that the image of $\varphi(C)$ is not a point.
\end{proof}

Let us show by means of an example that the condition in \Cref{lemma: relative positivity of -b_n on branch-chow over s_2 locus} stating that the generic point of $C$ lands in $\Branch^P_{\mathrm{Dn}}(\bP^N)$ cannot be removed.
\begin{ex} \label{example: CM not positive away from S2 locus}
    Set $n=N=3$, and consider the projective space $\bP^3$ with homogeneous coordinates $[x_0: x_1: x_2 : x_3]$. Let $H \subset \bP^3$ denote the linear subspace given by the vanishing of $x_3$, and let $\Sigma \subset H$ denote the line given by the vanishing of $x_2$ and $x_3$. Consider the family of closed subschemes $Z \subset \bP^3 \times \bP^1$ cut out by the three bihomogeneous equations $x_3^2$,  $x_3x_2$ and $x_3(x_1u - x_0v))$, where we denote the homogeneous coordinates of the $\mathbb{P}^1$ component by $[u:v]$. This is a $\mathbb{P}^1$-flat family of subschemes of $\bP^3$. The fiber $Z_p$ over a given point $p = [a:b]$ in $\bP^1$ is given by a subscheme of $\bP^3$ consisting of $H$ with an embedded point at $[x_0:x_1] = [a:b]$ in $\Sigma \subset H$. Consider the $\mathbb{P}^1$-family of reduced schemes $X:= (\bP^3 \times \bP^1) \cup_Z (\bP^3 \times \bP^1)$ obtained by gluing two copies of $\bP^3 \times \bP^1$ along $Z$. Projection onto the first coordinate $\bP^3$ induces a well-defined morphism $f:X \to \bP^3$, which is a nontrivial $\bP^1$-family of equidimensional branchvarieties over $\bP^3$. We thus obtained a morphism $\varphi: \bP^1 \to \Branch^3(\bP^3)$. Since $n=N=3$, the Branch-Chow morphism $\Chow: \Branch^3(\bP^3) \to \Spec(k)$ has as target a point, so this family lies over the unique fiber. Let $q \in \bP^3(k)$ be a point that does not lie in $H$. The fiber product $X \times_{\bP^3} q$ is isomorphic to a trivial $\bP^1$-family $\pi: \bP^1 \sqcup \bP^1 \to \bP^1$. By definition of $b_3$, we have that $\varphi^*(b_3) = \det(\pi_*(X \times_{\bP^3} q)) = \cO_{\bP^1}$ for this nontrivial morphism $\varphi: \bP^1 \to \Branch^3(\bP^3)$.
\end{ex}

\begin{thm} \label{thm: ampleness of cm line bundle on s_2 locus}
    Let $X$ be a projective $k$-scheme equipped with a very ample line bundle $\cO_X(1)$ which induces an embedding $X \hookrightarrow \mathbb{P}^N$.  Fix a Hilbert polynomial $P$ of degree $n$. Let $a>n$ be a rational number. Then the following hold:
    \begin{enumerate}
        \item  Let $f: \Branch^P(X) \to \mathbb{P}^J$ denote any morphism induced by a base-point free linear system of a power of the semiample $\mathbb{Q}$-line bundle $ab_{n+1} - b_n$. Then the restriction of $f$ to $\Branch_{\mathrm{Dn}}^P(X) \subset \Branch^P(X)$ is quasi-finite.
        \item The $\mathbb{Q}$-line bundle $ab_{n+1} - b_n$ is ample on $\Branch_{\mathrm{Dn}}^P(X)$.
         \item The line bundle $-b_n$ is relatively ample for $\Chow: \Branch^P_{\mathrm{Dn}}(X) \to \bP(V_{P,N})$.
    \end{enumerate}
\end{thm}
\begin{proof} All of the statements may be checked after base-changing to the algebraic closure of $k$, and hence we may assume without loss of generality that $k$ is algebraically closed.\\

\noindent \textit{Part (1).} Consider the projective coarse space $\Branch^P(X) \to M$ with open subscheme $M' \subset M$ corresponding to the deminormal locus. The morphism $\Branch^P(X) \to \mathbb{P}^J$ factors through a morphism $M \to \mathbb{P}^J$. Suppose for the sake of contradiction that the composition $\psi: M' \hookrightarrow M \to \mathbb{P}^J$ is not quasi-finite. Then, there is a (possibly non-proper) curve $C^{\circ} \hookrightarrow M'$ which is contracted by $\psi$. The closure $C \subset M$ of $C^{\circ}$ in $M$ is a projective curve that is contracted by the morphism $M \to \mathbb{P}^J$, and so the pullback of $\cO_{\mathbb{P}^J}(1)$ to $C$ is trivial. Choose a finite smooth cover $\widetilde{C} \to C$ such that the composition $\widetilde{C} \to C \to M$ lifts to a morphism $\widetilde{C} \to \Branch^P(X)$. Consider the composition $\varphi: \widetilde{C} \to \Branch^P(X) \hookrightarrow \Branch^P(\mathbb{P}^N)$. By construction $\deg(\varphi^*(ab_{n+1}-b_n)) =0$. However, since $\varphi(\widetilde{C}) \subset |\Branch^P(\mathbb{P}^N)|$ is not a single point and the image of the generic point of $\widetilde{C}$ is contained in the deminormal locus $\Branch_{\mathrm{Dn}}^P(\mathbb{P}^N) \subset \Branch^P(\mathbb{P}^N)$,  \Cref{lemma: relative positivity of -b_n on branch-chow over s_2 locus} and \Cref{prop: -bd nef on branch} jointly imply that the pullback of $\varphi^*(ab_{n+1} - b_n)$ is positive\footnote{We may write $ab_{n+1} - b_n = \frac{1}{2}(a-n)b_{n+1} + \left((n + \frac{1}{2}(a-n))b_{n+1}- b_n \right)$, where $\frac{1}{2}(a-n)b_{n+1}$ is the pullback of an ample line bundle via the Chow morphism, and $(n + \frac{1}{2}(a-n))b_{n+1}- b_n$ is nef. We have two possible cases for the curve $C \to M$. If the curve $C$ is not contracted by the Chow morphism, then the term $\frac{1}{2}(a-n)b_{n+1}$ is positive, while the second term is non-negative, so that the sum is positive. Otherwise, the curve must be contained in a fiber of the Chow morphism. In this case the pull-back of $b_{n+1}$ and $-b_n$ are trivial and ample, respectively, so the pull-back of $a b_{n+1} - b_n$ is positive.}, a contradiction.
\\

\noindent \textit{Part (2).} This follows from Part (1) by applying Zariski's main theorem to the quasi-finite morphism $\psi:M' \to \bP^J$ induced by a power of $ab_{n+1} - b_n$.\\

\noindent \textit{Part (3).} Since $b_{n+1}$ is the pullback $\Chow^*(\cO(1))$ under the Chow morphism $\Chow: \Branch^P(X) \hookrightarrow \Branch^P(\mathbb{P}^N) \to \bP(V_{P,N})$ (\Cref{prop: b_{n+1} is the pullback of O(1) under Chow}), the relative ampleness of $-b_n$ is equivalent to the relative ampleness of $ab_{n+1} - b_n$, which was proven in Part (2).
\end{proof}

\begin{rem}
    The nontrivial morphism $\varphi: \mathbb{P}^1 \to \Branch^3(\mathbb{P}^3)$ constructed in \Cref{example: CM not positive away from S2 locus} satisfies $\varphi^*(b_n) = \varphi^*(b_{n+1}) = \cO_{\mathbb{P}^1}$. In particular, $\varphi^*(ab_{n+1} - b_n)$ is trivial for all $a$. This shows that restricting to the demi-normal locus is necessary to ensure the ampleness of $ab_{n+1} - b_n$ in \Cref{thm: ampleness of cm line bundle on s_2 locus}(2).
\end{rem}

In the following example, we show that the hypothesis $a>n$ is necessary for \Cref{thm: ampleness of cm line bundle on s_2 locus}.
\begin{ex}
   Consider the Segre embedding $\iota: \mathbb{P}^1 \times \mathbb{P}^1 \hookrightarrow \mathbb{P}^3$ induced by the complete linear series of the line bundle $\cO(1) \boxtimes \cO(1)$ on $\mathbb{P}^1 \times \mathbb{P}^1$. The second projection $p_2: \mathbb{P}^1 \times \mathbb{P}^1 \to \mathbb{P}^1$ exhibits $\mathbb{P}^1 \times \mathbb{P}^1$ as a smooth family of varieties over $\mathbb{P}^1$. We may view the induced morphism $\iota \times p_2: \mathbb{P}^1 \times \mathbb{P}^1 \to \mathbb{P}^3 \times \mathbb{P}^1$ as a $\mathbb{P}^1$-family of branchvarieties, which corresponds to a morphism $\varphi: \mathbb{P}^1 \to \Branch^1_{\mathrm{Dn}}(\mathbb{P}^3)$. In this case $\varphi^*(\lambda_t) = \cO(t^2 +t)$, and hence we have $b_2 = b_1 = \cO(2)$. Therefore, $\varphi^*(ab_2 - b_1)$ is ample only if $a>1$.
\end{ex}

\section{Projective moduli spaces of ample linear series}

Our results in this paper make it possible to analyze the following moduli problem.

\begin{defn}[Stack of ample linear series]
   We define the stack $\LS_P^r$ of \textbf{ample linear series} of rank $r \geq 0$ and Hilbert polynomial $P(t)$ to be the pseudofunctor that sends a $k$-scheme $T$ to the groupoid of tuples $(X \to T, \cO_X(1), \cE, s)$, where
    \begin{enumerate}
        \item $\pi: X \to T$ is a flat projective morphism with reduced equidimensional fibers, equipped with a relatively ample line bundle $\cO_X(1)$ with Hilbert polynomial $P(t)$.
    
        \item $\cE$ is a rank $r+1$ vector bundle on $T$.

        \item $s: \pi^*(\cE) \twoheadrightarrow \cO_X(1)$ is a surjection of $\cO_{X}$-modules.
    \end{enumerate}
\end{defn}

Note that a geometric point of $\LS_P^r$ amounts to an equidimensional variety $X$ equipped with a finite morphism $f:X \to \mathbb{P}^r$ such that the Hilbert polynomial with respect to the polarization $\cO_X(1) := f^*(\cO_{\mathbb{P}^n}(1))$ is $P(t)$. In fact, we have $\LS_P^r = \Branch^P(\mathbb{P}^r)/\GL_{r+1}$, where the $\GL_{r+1}$-action is induced from its standard linear action on $\mathbb{P}^r$.

We may define the line bundles $\lambda_m$ and $b_i$ on the stack $\LS_P^r$ similarly as for $\Branch^P(\mathbb{P}^r)$. Let $n$ be the degree of $P(t)$, and write $P(t) = \sum_{i=0}^n a_i \binom{t}{i}$. Fix once and for all a constant $A>n$. (For definiteness, one may choose $A=n+1$.)

For any $s,t \in \bQ$, let $L(s,t)$ denote the $\bQ$-line bundle on $\LS_P^r$,
\begin{equation}\label{E:define_L(s,t)}
L(s,t) := s \left( A b_{n+1} - b_n - \frac{B}{r+1} \det(\cE) \right) + \left(\frac{1}{a_n} b_{n+1} - \frac{n+1}{t P(t)} \lambda_t \right),
\end{equation}
where $B := A (n+1) a_n - n (a_n+a_{n-1})$ is a constant, which has been chosen so that family of $\bQ$-line bundles $L(s,t)$ on $\LS_P^r = \Branch^P(\mathbb{P}^r)/\mathrm{GL}_{r+1}$ has weight zero with respect to the central copy of $\mathbb{G}_m$ inside $\mathrm{GL}_{r+1}$.

Using the line bundles $L(1,t)$, we define a notion of semistability for ample linear series (here we employ the notion of $\Theta$-semistability developed in \cite{instability}). For this definition, we first need to recall the following.
\begin{defn}[Test configuration] \label{defn: test configuration}
    Let $K \supset k$ be an algebraically closed field extension, and let $p: \Spec(K) \to \LS_P^r$ be a geometric point, corresponding to a tuple $(X, \cO_X(1), \cE, s)$ defined over $K$. A \textbf{test configuration of $p$} is a morphism of stacks $f: \mathbb{A}^1_K/\mathbb{G}_m \to  \LS_P^r$ with an identification $f(1) \cong p$. Equivalently, a test configuration of $p$ is a $\mathbb{G}_m$-equivariant $\mathbb{A}^1_K$-family $(\widetilde{X}, \cO_{\widetilde{X}}(1), \widetilde{\cE}, \widetilde{s})$ along with an identification of the $1$-fiber with $(X, \cO_X(1), \cE, s)$.
\end{defn}
Given a test configuration $f: \mathbb{A}^1_K/\mathbb{G}_m \to \LS_P^r$ of a geometric point $p \in |\LS_P^r|$, we denote by $\wt(f^*L(1,t))$ the weight of the $0$-fiber of the $\mathbb{G}_m$-equivariant $\bQ$-line bundle $f^*L(1,t)$ on $\mathbb{A}^1_K$.

\begin{defn}[Polynomial semistability]\label{D:polynomial_stability}
    We say that $p \in |\LS_P^r|$ is \textbf{polynomial semistable} if for all test configurations $f$ of $p$, we have $\wt(f^*L(1,t)) \geq 0$ for all sufficiently large $t \gg 0$.
\end{defn}

The following shows that points in $\LS_P^{r,\rm{ss}}$ correspond to linear series in the traditional sense, i.e., rank $r+1$ subspaces of $\Gamma(X,\cO_X(1))$.

\begin{lem}
    If $(X, \cO_X \otimes V \to \cO_X(1))$ corresponds to a semistable point of $\LS_P^r$, then the homomorphism $V \to \Gamma(X,\cO_X(1))$ is injective. In particular, the stack $\LS_P^{r,\rm{ss}}$ is representable over the stack of polarized reduced projective schemes.
\end{lem}
\begin{proof}
    Let $K = \ker(V \to \Gamma(X,\cO_X(1)))$. If $K \neq 0$, then consider the following filtration of the point in $\LS_P^r$: the equivariant degeneration of $(X,\cO_X(1))$ over $\bA^1 / \bG_m$ is trivial, and the filtration of $V$ assigns $F_w(V) = K$ for $w\geq 1$ and $F_w(V)=V$ otherwise. For this morphism $\bA^1 / \bG_m \to \LS_P^r$, $\lambda_t$, $b_n$, and $b_{n+1}$ all have weight $0$, because they are determined solely by the family of polarized varieties, and $\det(\cE)$ has weight $\dim(K)>0$. It follows that $\wt(L(1,t))<0$, so this filtration is destabilizing.
    
    The kernel of the homomorphism from the automorphism group of a point $(X,\cO_X \otimes V \to \cO_X(1)) \in \LS_P^r$ to the automorphism group of the underlying polarized variety $(X,\cO_X(1))$ consists of automorphisms $g \in \GL(V)$ that commute with the map $V \to \Gamma(X,\cO_X(1))$. If $V \to \Gamma(X,\cO_X(1))$ is injective, then this kernel is trivial, which implies that the forgetful functor from $\LS_P^r$ to the stack of polarized reduced projective schemes is representable.
\end{proof}

\begin{thm}\label{T:moduli_ample_linear_series}
    The locus of polynomial semistable ample linear series is an open substack $\LS_P^{r,\mathrm{ss}} \subset \LS_P^r$ which admits a projective good moduli space $M_P^r$. Furthermore, for all sufficiently large $ t \gg 0$, there exists $c>0$ such that for all $ s>c$, the restriction $L(s,t)|_{\LS_{P}^{r,\mathrm{ss}}}$ descends to an ample $\bQ$-line bundle on $M_P^r$.
\end{thm}

\begin{proof}

The theorem is a direct application of the main theorem of GIT for semiprojective DM stacks, which we have spelled out in \Cref{appendix}.

The second term of $L(s,t)$ converges in the real N\'eron-Severi group $\NS(\Branch^P(\bP^r))_\bR$ to $0$ as $t \to \infty$. This implies that the polynomial stability condition in \Cref{D:polynomial_stability} is lexicographic, i.e., a point is semistable if and only for every filtration, the $s$ term of \eqref{E:define_L(s,t)} has weight $\geq 0$, and if it has weight $0$, then the second term has weight $\geq 0$ for all $t \gg 0$. In particular, even though $s$ is set to $1$ in \Cref{D:polynomial_stability}, the resulting notion of semistability agrees with polynomial semistability as in \Cref{D:polynomial_semistability_general} with respect to the two-variable polynomial $\tilde{L}(t_1,t_2) = t_2 P(t_2) L(t_1,t_2)$. We include the factor of $t_2 P(t_2)$ to clear denominators and get a polynomial sequence of $\bQ$-line bundles.

\Cref{prop: semiampleness of a b_n+1 - b_n} implies that $A b_{n+1}-b_n$ is the pullback of an ample bundle under some morphism $\pi : \Branch^P(\bP^r) \to Y$, and \Cref{prop: -lambda_n is relatively ample} implies that $-\lambda_r$ is $\pi$-ample for all $r \gg 0$. It follows that for all sufficiently large $x_2 \gg 0$, there exists $c >0$ such that $\tilde{L}(x_1,x_2)$ restricts to an ample line bundle on $\Branch^P(\bP^r)$ for all $x_1>c$. In other words, $\tilde{L}(t_1,t_2)$ is asymptotically ample on $\Branch^P(\bP^r)$. The result now follows from \Cref{T:main_GIT}.
\end{proof}

\begin{ex}[Moduli of Noether-normalized varieties]
The minimal value of $r$ for which $\LS_P^r$ is nonempty is $r = n$. In this case, the stack $\LS_P^r$ parameterizes families of reduced projective schemes equipped with a Noether normalization $X \to \bP^n$, up to the action of $\GL_{n+1}$. The semistability condition on $\LS_P^r$ is closely related to $K$-semistability of the underlying variety $X$. Indeed, because $b_{n+1}$ is trivial on $\Branch^n(\bP^n)$, it must be equivariantly proportional to $\det(\cE)$ on each connected component. This implies that the coefficient of $s$ in \eqref{E:define_L(s,t)} has the form $q \cdot b_{n+1}-b_n$ for some $q\in \bQ$. The CM line bunlde $\CM$ is the only linear combination of $b_n$ and $b_{n+1}$ that is invariant under the action of $\bG_m$ by scaling the polarization. It follows that for some positive number $C>0$,
\[
L(1,t) = C \cdot \CM + \frac{1}{(n+1)a_n} b_{n+1} - \frac{1}{t P(t)} \lambda_t.
\]
The leading order term of this sequence is $C \cdot \CM$, so if $(X,\cO_X(1))$ is a $K$-stable variety satisfying Serre's condition $S_2$, the corresponding point of $\LS_P^r$ is semistable in our sense.\footnote{Here one needs to use the corrected definition of $K$-stability in \cite[Def. 2.4]{odaka-stability-II}, where we restrict to test configurations that are not almost trivial (i.e. there are not isomorphic to a product test configuration away from codimension 2). Note that, if $x$ is a point of $\LS_P^r$ whose underlying variety is $S_2$, then any test configuration of $x$ as in \Cref{defn: test configuration} will have $S_2$ total space, and hence it will be almost trivial if and only if it is a product test configuration.} For example, if $X$ is a KSB-stable variety of general type, then it is $K$-stable \cite[Thm. 1.1]{odaka-stability-I}, so for any choice of Noether normalization $X \to \bP^n$ induced by a power of the canonical divisor, the corresponding point of $\LS_P^n$ is semistable.
\end{ex}

\appendix

\section{Polynomial stability and GIT for stacks}
\label{appendix}

One of the first applications of GIT was to construct moduli spaces of semistable schemes as GIT quotients of Hilbert schemes of $\bP^N$. This is sometimes referred to as Hilbert stability. One of the interesting applications of the projectivity of $\Branch(\bP^N)$ is to introduce this as an analogous tool for studying semistability of varieties using GIT. However, because $\Branch(\bP^N)$ is a DM stack, we need a slight generalization of GIT for this purpose.

In this appendix, we develop the theory of GIT for actions of reductive groups $G$ on a semi-projective stack $X$, as in the following.

\begin{defn}
    A finite type algebraic stack $X$ over $k$ is \textbf{semi-projective} if it has finite inertia, admits a proper morphism $X \to \Spec(A)$ to a finite type affine $k$-scheme, and has a line bundle $L$ such that some power $L^m$ descends to an ample bundle on the coarse moduli space of $X$. Such an $L$ is called ample.
\end{defn} 

\begin{defn} \label{defn: group action on a stack}
   If $G$ is a group scheme locally of finite type over $k$, \textbf{a $G$-action} on an algebraic stack $X$ is an algebraic stack $\cX$ equipped with a morphism $\cX \to BG$ and an isomorphism $X \cong \Spec(k) \times_{BG} \cX$. In this case, we say $\cX \cong X/G$ is the quotient stack associated to the $G$-action. 
\end{defn}

 \begin{ex}
 Let $X$ be a $k$-scheme equipped with the action of an algebraic group $G$ over $k$. The canonical morphism of stacks $q: X \to X / G$ corresponds to the data of the trivial $G$-bundle $G \times_{\Spec(k)} X$ and the action map $G \times_{\Spec(k)} X \to X$. The morphism $q$ and the structure map $X \to \Spec(k)$ induce an isomorphism $X \cong \Spec(k) \times_{BG} X / G$. Thus \Cref{defn: group action on a stack} generalizes the usual notion of a group action.
 \end{ex}

We regard $\bQ[t_1,\ldots,t_N]$ as an ordered vector space by using the usual ordering on $\bQ$ if $N=0$ and declaring recursively that $f(t_1,\ldots,t_N) > 0$ if there exists $ c \in \bQ$ such that for all $x \geq c$, the polynomial $f(t_1,t_2,\ldots,t_{N-1},x) \in \bQ[t_1,\ldots,t_{N-1}]$ is $> 0$. This is equivalent to saying that the leading coefficient of $f$ with respect to the lexicographic monomial order on $\bQ[t_1,\ldots,t_N]$ is positive.

Another way to characterize this ordering is to declare that a statement holds for all sufficiently large $ (x_1,\ldots,x_N) \gg 0$ in lex order if there are functions $c_i(t_{i+1},\ldots,t_N) :\mathbb{Q}^{N-1} \to \mathbb{Q}$ for $i=1,\ldots,N$ such that the statement holds for all $(x_1,\ldots,x_N) \in \bQ^N$ with $x_i > c_i(x_{i+1},\ldots,x_N)$. Then the ordering above is equivalent to saying that for $f \in \bQ[t_1,\ldots,t_N]$, $f > 0$ if $f(x_1,\ldots,x_n)>0$ for all $(x_1,\ldots, x_N) \gg 0$ in lex order.

A \textbf{polynomial sequence} of $\bQ$-line bundles $\cX$ is an element $L(t_1,\ldots,t_N) \in \Pic(\cX) \otimes \bQ[t_1,\ldots,t_N]$ for some $N>0$. For any tuple $x_1,\ldots,x_N \in \bQ$, we can evaluate such a polynomial to obtain $L(x_1,\ldots,x_N) \in \Pic(\cX)_\bQ$. We define $L(t_1,\ldots,t_N) \in \Pic(\cX) \otimes \bQ[t_1,\ldots,t_N]$ to be \emph{asymptotically ample} if $\forall (x_1,\ldots,x_N) \gg 0$ in lex order, $L(x_1,\ldots,x_N)$ is ample. This notion of ampleness also has a recursive description and a description in terms of leading monomials analogous to the notion of positivity in $\bQ[t_1,\ldots,t_N]$.


\begin{defn}\label{D:polynomial_semistability_general}
    Given a polynomial sequence $L = L(t_1,\ldots,t_N) \in \Pic(\cX)\otimes \bQ[t_1,\ldots,t_N]$, we say that a point $p \in |\cX|$ is \textbf{$L$-semistable} if for any field extension $K \supset k$ and any morphism $f : \Theta_K \to \cX$ with $f(1) = p$, the weight $\wt(f^\ast(L(t_1,\ldots,t_N))) \in \bQ[t_1,\ldots,t_N]$ is $\geq 0$.
\end{defn}

In addition to identifying semistable points, one can ask if unstable points have ``optimally destabilizing filtrations.'' In classical GIT, this is the Kempf-Hesselink optimally destabilizing cocharacter of an unstable point \cite{kempf1978instability}. It leads to a $G$-equivariant stratification of the unstable locus by locally closed subschemes whose points have the same optimally destabilizing cocharacter, up to conjugacy.

The theory of $\Theta$-stability and $\Theta$-stratifications provides an intrinsic generalization of this story. If one fixes a real-valued Weyl-invariant rational quadratic norm on the coweight lattice of $G$, or equivalently a rational quadratic norm $\|\bullet\|$ on graded points of $BG$ \cite[Defn. 4.1.12]{instability}, then for any filtration $f : \Theta_K \to \cX$, one can define the $\bR[t_1,\ldots,t_N]$-valued numerical invariant
\begin{equation} \label{E:numerical_invariant}
    \mu(f) = \frac{\wt(f^\ast(L(t_\bullet)))}{\lVert f\rVert},
\end{equation}
where $\lVert f \rVert$ denotes the norm of the cocharacter of $G$ corresponding to the composition $\Theta_K \to \cX \to BG$.

A point $x \in |\cX|$ is unstable if and only if there is a filtration with $f(1) = x$ and $\mu(f) < 0$. A HN filtration of $x$ one that minimizes $\mu(f)$ subject to the constraint $f(1)=x$, and we say that $\mu$ determines a $\Theta$-stratification on $\cX$ if the locus of HN filtrations is open in $\mathop{Map}(\Theta,\cX)$, and the components of this open substack can be identified with locally closed strata in $\cX$ under the forgetful map $f \mapsto f(1)$. See \cite[Def.~2.2.1]{instability} for a more precise discussion.

\begin{rem}
    When $X$ is a scheme and $N=0$, so that $L \in \Pic(\cX) \otimes \bQ$, \Cref{D:polynomial_semistability_general} is an intrinsic reformulation of the Hilbert-Mumford criterion for semistability. The more general formulation of polynomial semistability is a convenient way to avoid an arbitrary choice of a line bundle in situations where one only has a canonical asymptotically ample polynomial sequence. As we will see below, using polynomial semistability does not change the semistable locus at all, but it does allow one to make canonical choices of optimally destabilizing cocharacters.
\end{rem}

\begin{thm}\label{T:main_GIT}
Let $X$ be a semi-projective algebraic stack over $k$ equipped with an action of a reductive group $G$, and let $L \in \Pic(\cX) \otimes \bQ[t_1,\ldots,t_N]$ be an asymptotically ample polynomial sequence of line bundles on $\cX = X/G$. Then:
\begin{enumerate}
\item The numerical invariant \eqref{E:numerical_invariant} determines a $\Theta$-stratification of $\cX$.\\
\item The set of $L$-semistable points $\cX^{\rm ss} \subset \cX$ is open and admits a good moduli space $M$ that is proper over $\Spec(\Gamma(\cX,\cO_\cX))$.\\
\item For all sufficiently large $(x_1,\ldots,x_N) \gg 0$ in lex order:\\ 
\begin{enumerate}
    \item The polynomial semistable locus for $L(t_1,\ldots,t_N)$ agrees with the semistable locus with respect to $L(x_1,\ldots,x_N)$;\smallskip
    \item $L(x_1,\ldots,x_N)$ descends to an ample $\bQ$-line bundle on $M$; and\smallskip
    \item If the leading coefficient $L_{\max}$ of $L(t_1,\ldots,t_N)$ is semiample on $X$, then it descends to a semiample $\bQ$-line bundle on $M$ as well.
\end{enumerate}
\end{enumerate}
\end{thm}

\begin{proof}
We first claim that the numerical invariant $\mu$ on $\cX$ satisfies the monotonicity conditions introduced in \cite[Def.~5.2.1, Def.~5.5.7]{instability}. Let $R$ be a discrete valuation ring, and let $\cT$ be either of the stacks $\Spec(R[s,t]/(st-\pi)) / \bG_m$ or $\Spec(R[t])/\bG_m$, where $\bG_m$ acts with weight $-1$ on $t$ and $1$ on $s$. $\cT$ has a unique closed codimension-$2$ point $0 \in \cT$. Consider the projective morphism $X/G \to \Spec(A)/G$ for which $L(t_\bullet)$ is asymptotically relatively ample. For any morphism $\cT \setminus 0 \to \cX$, the composition $\cT \setminus 0 \to \cX \to \Spec(A)/G$ extends uniquely to a morphism $\cT \to \Spec(A)/G$, because the latter is $\Theta$-reductive and $S$-complete \cite[Prop.3.21(2)+Prop.3.44(2)]{existence_gms}. Now the fiber product $\cY := \cX \times_{\Spec(A)/G} \cT \to \cT$ is projective over $\cT$, and the pullback of $L(t_\bullet)$ is asymptotically relatively ample on $\cY$. The original morphism $\cT \setminus 0 \to \cX$ over $\Spec(A)/G$ gives a section $\cT \setminus 0 \to \cY$. We let $\cW$ be the closure of this section, equipped with its projection $p : \cW \to \cT$, and let $\tilde{f} : \cW \to \cX$ be the restriction of the morphism $\cY \to \cX$ to $\cW$. 

Next, consider a commutative diagram
\[
\xymatrix{ \bP^1_K / \bG_m \ar[r]^\phi \ar[d] & \cW \ar[d]^p \ar^{\tilde{f}}[r] & \cX \\ (B\bG_m)_K \ar[r] & \cT & },
\]
where the lower horizontal arrow is a positive multiple of the canonical graded point at $0 \in \cT$, and the induced map $\bP^1_K \to (B\bG_m)_K \times_{\cT} \cW$ is finite. By construction $\tilde{f}^\ast(L(t_\bullet))$ is asymptotically ample relative to $\cT$, so $\phi^\ast (\tilde{f}^\ast(L(t_\bullet)))$ is ample on $\bP^1_K$. This implies that the weight of this line bundle at $\infty$ is less that the weight at $0$. On the other hand, both graded points $\{0\} / \bG_m \to \cX$ and $\{\infty\}/\bG_m \to \cX$ have the same norm, because their compositions to $BG$ are both isomorphic to a positive multiple of the canonical graded point at $0 \in \cT$. So we have $\mu(\{\infty\}/\bG_m \to \cX) < \mu(\{0\} / \bG_m \to \cX)$, which is precisely the monotonicity condition \cite[Def.~5.2.1]{instability} (in the reference, the inequality goes the other way, because \cite{instability} uses the opposite sign convention, where filtrations with positive weight are destabilizing). 

We have now verified the hypotheses of the main theorem of the beyond geometric invariant theory program \cite[Thm.~5.5.10]{instability}: $\mu$ is $S$-monotone and $\Theta$-monotone, HN boundedness (condition (B)) holds because $\cX$ is quasi-compact, and the rest of the hypotheses follow because we are using a numerical invariant of the form \eqref{E:numerical_invariant}. Conclusions (1) and (2) of the theorem follow immediately from \cite[Thm.~5.5.10]{instability}. 

Next, consider the relative coarse moduli space morphism $\cX \to \cY$ over $BG$, where $\cY \cong Y / G$ and $Y$ is the coarse moduli space of $X$. Pullback along this morphism induces an isomorphism $\Pic(\cY)_\bQ \cong \Pic(\cX)_\bQ$,
\footnote{If $\pi : \cX \to \cY$ is a relative coarse moduli space morphism, the pullback functor on Picard groups $\pi^\ast : \Pic(\cY) \to \Pic(\cX)$ is injective. Indeed, $\pi$ is a good moduli space morphism since $k$ has characteristic $0$, so $F \cong \pi_{\ast}(\pi^{\ast}(F))$ for any $F \in \QCoh(\cX)$ by \cite[Prop. 4.5]{alper_good_moduli}. Therefore, $L \cong \pi_\ast(\pi^\ast(L))$ for any line bundle $L$ on $\cY$, and $\pi_\ast(\cO_\cX) \cong \cO_\cY$, so any line bundle that pulls back to $\cO_\cX$ was already isomorphic to $\cO_\cY$. Tensoring with $\bQ$ is exact, so we have an injection $\Pic(\cY)_\bQ \hookrightarrow \Pic(\cX)_\bQ$. This is surjective assuming $\cX$ is quasi-compact, because for any $L$ on $\cX$, there is some sufficiently divisible $m$ so that all automorphism groups act trivially on fibers of $L^{\otimes m}$, and this is precisely the condition for $L^{\otimes m} \cong \pi^\ast(L')$ for some $L'$ on $\cY$.}
so we can regard $L \in \Pic(\cY) \otimes \bQ[t_1,\ldots,t_N]$. Every filtration in $\cY$ lifts uniquely up to composition with a ramified covering $\Theta \to \Theta$ by \cite[Thm.~B.1]{gauged_maps}, so $\cX^{\rm ss}$ is the preimage of $\cY^{\rm ss}$ and the good moduli spaces of $\cX^{\rm ss}$ and $\cY^{\rm ss}$ are isomorphic. We may therefore assume for the remainder of the proof that $X$ is a semi-projective scheme. By embedding $G$ in $\GL_n$, we have $\cX \cong X' / \GL_n$, where $X' = (X \times \GL_n) / G$ is still semi-projective and the induced line bundle $L(t_1,\ldots,t_N)$ is still asymptotically ample. So in addition, we may assume $G = \GL_n$. Let $T \subset \GL_n$ be the standard maximal torus, and let $\mathbf{N}$ denote the coweight lattice of $T$.

Now every unstable point -- either with respect to the sequence $L(t_\bullet)$ or a specific ample line bundle $L(x_1,\ldots,x_N)$ -- lies in the image of $A_{Z,\lambda} / T \to \cX$, where $\lambda : \bG_m \to T$ is a non-trivial cocharacter with $\wt_\lambda(L|_Z)<0$, $Z \subset X^{\lambda}$ is a connected component of the $\lambda$-fixed points $X^{\lambda} \subset X$, and $A_{Z,\lambda} \hookrightarrow X$ is the locally closed subscheme of points $x$ for which $\lim_{z \to 0} \lambda(z) \cdot x \in Z$.
In fact, only finitely many pairs of subschemes $(Z,A)$ appear as $(Z,A_{Z,\lambda})$ for some cocharacter $\lambda$.\footnote{If one chooses an equivariant embedding $X \hookrightarrow \bP^n \times \bA^m$ with a linear action on the latter, then any pair $(Z,A_{Z,\lambda})$ is a connected component of $(X \cap Z', X \cap A'_{Z',\lambda})$, where $(Z',A'_{Z',\lambda})$ is a pair of fixed locus and attracting locus in $\bP^n \times \bA^m$. It therefore suffices to show this for $\bP(V) \times U$. In this case, both $Z'$ and $A'_{Z',\lambda}$ are linear subspaces spanned by a subset of the weight spaces of the linear representations $V$ and $U$. There are only finitely many spaces that arise as sums of weight spaces, so there are only finitely many pairs.} 
Let $I$ denote the set of such pairs.

For each $(Z,A) \in I$, let $\sigma_{Z,A} \subset \mathbf{N}_\bR$ be the closed rational polyhedral convex cone spanned by cocharacters $\lambda$ such that $Z \subset X^{\lambda}$ and $A$ is contained in an attracting stratum of $\lambda$. For any $(x_1,\ldots,x_N) \in \bQ^N$, the function $\wt_\lambda (L(x_\bullet)|_Z)$ defines a linear form $\ell_x : \Span(\sigma_{Z,A}) \to \bR$, and a cocharacter $\lambda \in \sigma_{Z,A}$ is destabilizing for points in the image of $A/T \to \cX$ with respect to the line bundle $L(x_1,\ldots,x_N)$ if and only if $\ell_x(\lambda)<0$. By linearity of $\ell_x$, this holds for some $\lambda \in \sigma_{Z,A}$ if and only if it holds for some $\lambda$ lying on an extremal ray of $\sigma_{Z,A}$.


Likewise, consider the polynomial-valued linear form $\ell_t: \Span(\sigma_{A,Z}) \to \bR[t_1,\ldots,t_N]$ given by
\[
\ell_t(\lambda) := \wt_{\lambda}(L(t_\bullet)|_Z) = \sum_\alpha f_\alpha(\lambda) t_1^{\alpha_1} \cdots t_N^{\alpha_N}
\]
where the coefficients $f_\alpha$ are rational linear forms on $\Span(\sigma_{Z,A})$. 
Because the positive cone in $\bQ[t_1,\ldots,t_N]$ is closed under linear combinations with positive real coefficients, the same reasoning as in the previous paragraph shows that there is some $\lambda \in \sigma_{Z,A}$ such that $\ell_t(\lambda) < 0$ in $\bQ[t_1,\ldots,t_N]$ if and only if there is an extremal vector of $\sigma_{Z,A}$ with this property.

We have shown that to check whether $\sigma_{Z,A}$ contains a destabilizing $\lambda$ for either $\ell_x$ or $\ell_t$, it suffices to test a finite set of $\lambda$ that generate all extremal rays of $\sigma_{Z,A}$. For any fixed $\lambda$, however, we have that for all $(x_1,\ldots,x_N) \gg 0$ in lex order, the sign of $\sum_\alpha f_\alpha(\lambda) x_1^{\alpha_1} \cdots x_N^{\alpha_N}$ agrees with the sign of $f_{\alpha^*}(\lambda)$, where $\alpha^*$ is the maximal $\alpha$ in lexicographic monomial order for which the coefficient $f_{\alpha}(\lambda) \neq 0$. This in turn agrees with the sign of $\sum_\alpha f_\alpha(\lambda) t_1^{\alpha_1} \cdots t_N^{\alpha_N}$ in $\bQ[t_1,\ldots,t_N]$. We conclude that for all $ (x_1,\ldots,x_N) \gg 0$ in lex order, $\sigma_{Z,A}$ contains a destabilizing cocharacter for $\ell_x$ if and only if it contains a destabilizing cocharacter for $\ell_t$. Because there are only finitely many $(Z,A) \in I$, this shows that the unstable locus for the polynomial $L(t_1,\ldots,t_n)$ agrees with the unstable locus for $L(x_1,\ldots,x_N)$ for all $(x_1,\ldots,x_N) \gg 0$ in lex order.

Having proved 3(a), item 3(b) is the standard formulation of the main theorem of GIT \cite[Thm. 10]{mumford1994geometric} using the Hilbert-Mumford criterion \cite[Thm. 2.1]{mumford1994geometric}. As long as $m$ is a multiple of $\pi_0(G_x)$ for every stabilizer group $G_x$ of a closed orbit in $X^{\rm ss}$, $L_r^m$ descends to the good moduli space, and it is ample there.

To prove 3(c), suppose $L_\alpha \in \Pic(\cX)$ is the leading coefficient of $L(t_1,\ldots,t_N)$. Let $Y = \Proj(\bigoplus_{i \geq 0} \Gamma(X,L_\alpha^{\otimes i}))$. If $L_\alpha$ is semiample, there is a projective morphism $X \to Y$ such that some multiple of $L_\alpha$ is the pullback of an ample line bundle on $Y$. A point in $X$ has a limit in $X$ under a cocharacter $\lambda : \bG_m \to G$ if and only if its image in $Y$ has a limit. It follows that 
the polynomial semistable locus $X^{\rm ss} \subset X$ is contained in the preimage under $X \to Y$ of the $L_\alpha$-semistable locus $Y^{\rm{ss}} \subset Y$. One therefore obtains a morphism of good moduli spaces $X^{\rm ss} /\!/ G \to Y^{\rm ss} /\!/ G$. It follows that $m L_\alpha$, which descends to the GIT quotient $X^{\rm ss} /\!/ G$, is pulled back from an ample bundle on $Y^{\rm ss}/\!/G$, hence it is semiample on $X^{\rm ss} /\!/ G$.
\end{proof}


\footnotesize
\bibliography{moduli_algebras}

\begin{bibdiv}
\begin{biblist}

\bib{alper_good_moduli}{article}{
      author={Alper, Jarod},
       title={Good moduli spaces for {A}rtin stacks},
        date={2013},
        ISSN={0373-0956,1777-5310},
     journal={Ann. Inst. Fourier (Grenoble)},
      volume={63},
      number={6},
       pages={2349\ndash 2402},
         url={https://doi.org/10.5802/aif.2833},
      review={\MR{3237451}},
}

\bib{existence_gms}{article}{
      author={Alper, Jarod},
      author={Halpern-Leistner, Daniel},
      author={Heinloth, Jochen},
       title={Existence of moduli spaces for algebraic stacks},
        date={2023},
        ISSN={0020-9910,1432-1297},
     journal={Invent. Math.},
      volume={234},
      number={3},
       pages={949\ndash 1038},
         url={https://doi.org/10.1007/s00222-023-01214-4},
      review={\MR{4665776}},
}

\bib{branchvarieties-paper}{article}{
      author={Alexeev, Valery},
      author={Knutson, Allen},
       title={Complete moduli spaces of branchvarieties},
        date={2010},
        ISSN={0075-4102,1435-5345},
     journal={J. Reine Angew. Math.},
      volume={639},
       pages={39\ndash 71},
         url={https://doi.org/10.1515/CRELLE.2010.010},
      review={\MR{2608190}},
}

\bib{stble_maps_git}{article}{
      author={Baldwin, Elizabeth},
      author={Swinarski, David},
       title={A geometric invariant theory construction of moduli spaces of stable maps},
        date={2008},
        ISSN={1687-3017,1687-3009},
     journal={Int. Math. Res. Pap. IMRP},
      number={1},
       pages={Art. ID rpn 004, 104},
      review={\MR{2431236}},
}

\bib{codogni2021positivity}{article}{
      author={Codogni, Giulio},
      author={Patakfalvi, Zsolt},
       title={Positivity of the cm line bundle for families of k-stable klt fano varieties},
        date={2021},
     journal={Inventiones mathematicae},
      volume={223},
       pages={811\ndash 894},
}

\bib{eisenbud-harris}{article}{
      author={Eisenbud, David},
      author={Harris, Joe},
       title={Limit linear series: basic theory},
        date={1986},
        ISSN={0020-9910,1432-1297},
     journal={Invent. Math.},
      volume={85},
      number={2},
       pages={337\ndash 371},
         url={https://doi.org/10.1007/BF01389094},
      review={\MR{846932}},
}

\bib{fine2006note}{article}{
      author={Fine, Joel},
      author={Ross, Julius},
       title={A note on positivity of the cm line bundle},
        date={2006},
     journal={International Mathematics Research Notices},
      volume={2006},
       pages={O95875},
}

\bib{rho-sheaves}{article}{
      author={G\'omez, Tom\'as~L.},
      author={Fernandez~Herrero, Andres},
      author={Zamora, Alfonso},
       title={The moduli stack of principal {$\rho$}-sheaves and {G}ieseker-{H}arder-{N}arasimhan filtrations},
        date={2024},
        ISSN={0025-5874,1432-1823},
     journal={Math. Z.},
      volume={307},
      number={3},
       pages={Paper No. 51, 67},
         url={https://doi.org/10.1007/s00209-024-03497-6},
      review={\MR{4756614}},
}

\bib{gieseker_moduli_curves}{book}{
      author={Gieseker, D.},
       title={Lectures on moduli of curves},
      series={Tata Institute of Fundamental Research Lectures on Mathematics and Physics},
   publisher={Tata Institute of Fundamental Research, Bombay; by Springer-Verlag, Berlin-New York},
        date={1982},
      volume={69},
        ISBN={3-540-11953-1},
      review={\MR{691308}},
}

\bib{fga_grothendieck}{incollection}{
      author={Grothendieck, Alexander},
       title={Fondements de la g\'eom\'etrie alg\'ebrique. {Commentaires}},
    language={fr},
        date={1962},
   booktitle={S\'eminaire bourbaki : ann\'ee 1961/62, expos\'es 223-240},
      series={S\'eminaire Bourbaki},
   publisher={Soci\'et\'e math\'ematique de France},
       pages={297\ndash 298},
         url={http://www.numdam.org/item/SB_1961-1962__7__297_0/},
}

\bib{halpern2024moduli}{article}{
      author={Halpern-Leistner, Daniel},
      author={Fernandez~Herrero, Andres},
      author={Jones, Trevor},
       title={Moduli spaces of sheaves via affine grassmannians},
        date={2024},
     journal={Journal f{\"u}r die reine und angewandte Mathematik (Crelles Journal)},
      number={0},
}

\bib{instability}{misc}{
      author={Halpern-Leistner, Daniel},
       title={On the structure of instability in moduli theory},
         how={\url{https://arxiv.org/abs/1411.0627}},
        date={2022},
}

\bib{gauged_maps}{misc}{
      author={Halpern-Leistner, Daniel},
      author={Herrero, Andres~Fernandez},
       title={The structure of the moduli of gauged maps from a smooth curve},
        date={2023},
         url={https://arxiv.org/abs/2305.09632},
}

\bib{huybrechts2010geometry}{book}{
      author={Huybrechts, Daniel},
      author={Lehn, Manfred},
       title={The geometry of moduli spaces of sheaves},
   publisher={Cambridge University Press},
        date={2010},
}

\bib{kempf1978instability}{article}{
      author={Kempf, George~R},
       title={Instability in invariant theory},
        date={1978},
     journal={Annals of Mathematics},
      volume={108},
      number={2},
       pages={299\ndash 316},
}

\bib{kollar_projectivity}{article}{
      author={Koll\'ar, J\'anos},
       title={Projectivity of complete moduli},
        date={1990},
        ISSN={0022-040X,1945-743X},
     journal={J. Differential Geom.},
      volume={32},
      number={1},
       pages={235\ndash 268},
         url={http://projecteuclid.org/euclid.jdg/1214445046},
      review={\MR{1064874}},
}

\bib{kollar-singularities-mmp}{book}{
      author={Koll\'ar, J\'anos},
       title={Singularities of the minimal model program},
      series={Cambridge Tracts in Mathematics},
   publisher={Cambridge University Press, Cambridge},
        date={2013},
      volume={200},
        ISBN={978-1-107-03534-8},
         url={https://doi.org/10.1017/CBO9781139547895},
        note={With a collaboration of S\'andor Kov\'acs},
      review={\MR{3057950}},
}

\bib{kollar-families-varieties-book}{book}{
      author={Koll\'ar, J\'anos},
       title={Families of varieties of general type},
      series={Cambridge Tracts in Mathematics},
   publisher={Cambridge University Press, Cambridge},
        date={2023},
      volume={231},
        ISBN={978-1-009-34610-8},
        note={With the collaboration of Klaus Altmann and S\'andor J. Kov\'acs},
      review={\MR{4566297}},
}

\bib{knudsen-mumford}{article}{
      author={Knudsen, Finn~Faye},
      author={Mumford, David},
       title={The projectivity of the moduli space of stable curves. {I}. {P}reliminaries on ``det'' and ``{D}iv''},
        date={1976},
        ISSN={0025-5521},
     journal={Math. Scand.},
      volume={39},
      number={1},
       pages={19\ndash 55},
         url={https://doi-org.proxy.library.cornell.edu/10.7146/math.scand.a-11642},
      review={\MR{437541}},
}

\bib{liu-xu-zhuang}{article}{
      author={Liu, Yuchen},
      author={Xu, Chenyang},
      author={Zhuang, Ziquan},
       title={Finite generation for valuations computing stability thresholds and applications to {K}-stability},
        date={2022},
        ISSN={0003-486X,1939-8980},
     journal={Ann. of Math. (2)},
      volume={196},
      number={2},
       pages={507\ndash 566},
         url={https://doi.org/10.4007/annals.2022.196.2.2},
      review={\MR{4445441}},
}

\bib{mumford1994geometric}{book}{
      author={Mumford, David},
      author={Fogarty, John},
      author={Kirwan, Frances},
       title={Geometric invariant theory},
   publisher={Springer Science \& Business Media},
        date={1994},
      volume={34},
}

\bib{miranda-elliptic-p1}{article}{
      author={Miranda, Rick},
       title={The moduli of {W}eierstrass fibrations over {${\bf P}\sp{1}$}},
        date={1981},
        ISSN={0025-5831,1432-1807},
     journal={Math. Ann.},
      volume={255},
      number={3},
       pages={379\ndash 394},
         url={https://doi.org/10.1007/BF01450711},
      review={\MR{615858}},
}

\bib{miranda-elliptic-lectures}{book}{
      author={Miranda, Rick},
       title={The basic theory of elliptic surfaces},
      series={Dottorato di Ricerca in Matematica. [Doctorate in Mathematical Research]},
   publisher={ETS Editrice, Pisa},
        date={1989},
      review={\MR{1078016}},
}

\bib{mumford-stability-projective-varieties}{article}{
      author={Mumford, David},
       title={Stability of projective varieties},
        date={1977},
        ISSN={0013-8584},
     journal={Enseign. Math. (2)},
      volume={23},
      number={1-2},
       pages={39\ndash 110},
      review={\MR{450272}},
}

\bib{odaka-stability-I}{article}{
      author={Odaka, Yuji},
       title={The {C}alabi conjecture and {K}-stability},
        date={2012},
        ISSN={1073-7928,1687-0247},
     journal={Int. Math. Res. Not. IMRN},
      number={10},
       pages={2272\ndash 2288},
         url={https://doi.org/10.1093/imrn/rnr107},
      review={\MR{2923166}},
}

\bib{odaka-stability-II}{article}{
      author={Odaka, Yuji},
       title={The {GIT} stability of polarized varieties via discrepancy},
        date={2013},
        ISSN={0003-486X,1939-8980},
     journal={Ann. of Math. (2)},
      volume={177},
      number={2},
       pages={645\ndash 661},
         url={https://doi.org/10.4007/annals.2013.177.2.6},
      review={\MR{3010808}},
}

\bib{osserman-linear-series}{article}{
      author={Osserman, Brian},
       title={A limit linear series moduli scheme},
        date={2006},
        ISSN={0373-0956,1777-5310},
     journal={Ann. Inst. Fourier (Grenoble)},
      volume={56},
      number={4},
       pages={1165\ndash 1205},
         url={https://doi.org/10.5802/aif.2209},
      review={\MR{2266887}},
}

\bib{rydh2003chow}{thesis}{
      author={Rydh, David},
       title={Chow varieties},
        type={Master's Thesis},
        date={2003},
        note={Available at: \url{https://people.kth.se/~dary/Chow.pdf}},
}

\bib{rydh2008families}{misc}{
      author={Rydh, David},
       title={Families of cycles},
        date={2008},
        note={Part of thesis available at: \url{https://people.kth.se/~dary/famofcycles20080518.pdf}},
}

\bib{rydh-approximation}{article}{
      author={Rydh, David},
       title={Noetherian approximation of algebraic spaces and stacks},
        date={2015},
        ISSN={0021-8693,1090-266X},
     journal={J. Algebra},
      volume={422},
       pages={105\ndash 147},
         url={https://doi.org/10.1016/j.jalgebra.2014.09.012},
      review={\MR{3272071}},
}

\bib{ross-thomas}{article}{
      author={Ross, Julius},
      author={Thomas, Richard},
       title={A study of the {H}ilbert-{M}umford criterion for the stability of projective varieties},
        date={2007},
        ISSN={1056-3911,1534-7486},
     journal={J. Algebraic Geom.},
      volume={16},
      number={2},
       pages={201\ndash 255},
         url={https://doi.org/10.1090/S1056-3911-06-00461-9},
      review={\MR{2274514}},
}

\bib{mo186650}{misc}{
      author={Schwede, Karl},
       title={Obtaining non-normal varieties by pushout},
        date={2014},
        note={MathOverflow answer available at: \url{https://mathoverflow.net/questions/186406/obtaining-non-normal-varieties-by-pushout}},
}

\bib{SGA1}{misc}{
      author={{SGA1}},
       title={Rev\^etements \'etales et groupe fondamental},
      series={Documents Math\'ematiques (Paris)},
   publisher={Soci\'et\'e{} Math\'ematique de France, Paris},
        date={2003},
      volume={3},
        note={S\'eminaire de g\'eom\'etrie alg\'ebrique du Bois Marie 1960--61. Directed by A. Grothendieck, with two papers by M. Raynaud. Updated and annotated reprint of the 1971 original},
}

\bib{stacks-project}{misc}{
      author={{Stacks Project Authors}, The},
       title={\textit{Stacks Project}},
         how={\url{https://stacks.math.columbia.edu}},
        date={2024},
}

\bib{schurg2015derived}{article}{
      author={Sch{\"u}rg, Timo},
      author={To{\"e}n, Bertrand},
      author={Vezzosi, Gabriele},
       title={Derived algebraic geometry, determinants of perfect complexes, and applications to obstruction theories for maps and complexes},
        date={2015},
     journal={Journal f{\"u}r die reine und angewandte Mathematik (Crelles Journal)},
      volume={2015},
      number={702},
       pages={1\ndash 40},
}

\bib{bertini-theorems-tanaka}{article}{
      author={Tanaka, Hiromu},
       title={Bertini theorems admitting base changes},
        date={2024},
        ISSN={0021-8693,1090-266X},
     journal={J. Algebra},
      volume={644},
       pages={64\ndash 125},
         url={https://doi.org/10.1016/j.jalgebra.2023.12.038},
}

\bib{resolution_singularities}{article}{
      author={W{\l}odarczyk, Jaros{\l}aw},
       title={Simple {H}ironaka resolution in characteristic zero},
        date={2005},
        ISSN={0894-0347,1088-6834},
     journal={J. Amer. Math. Soc.},
      volume={18},
      number={4},
       pages={779\ndash 822},
         url={https://doi.org/10.1090/S0894-0347-05-00493-5},
      review={\MR{2163383}},
}

\bib{xu-fano-book}{misc}{
      author={Xu, Chenyang},
       title={K-stability of {F}ano varieties},
        date={2024},
        note={Unpublished book available at: \url{https://web.math.princeton.edu/~chenyang/Kstabilitybook.pdf}},
}

\bib{chenyang_positive}{article}{
      author={Xu, Chenyang},
      author={Zhuang, Ziquan},
       title={On positivity of the {CM} line bundle on {K}-moduli spaces},
        date={2020},
        ISSN={0003-486X,1939-8980},
     journal={Ann. of Math. (2)},
      volume={192},
      number={3},
       pages={1005\ndash 1068},
         url={https://doi.org/10.4007/annals.2020.192.3.7},
      review={\MR{4172625}},
}

\end{biblist}
\end{bibdiv}
\bibliographystyle{alpha}

\end{document}